\documentclass[10pt,a4paper]{amsart}

\usepackage{amsmath}
\usepackage{amssymb}
\usepackage{amsthm}
\usepackage{amsfonts}
\usepackage{microtype}
\usepackage{mathrsfs}
\usepackage{graphicx}
\usepackage{color}
\usepackage{latexsym}
\usepackage{rotating}
\usepackage{xspace}
\usepackage[all]{xy}
\usepackage{longtable}
\usepackage{leftidx}
\usepackage{mathtools}
\usepackage{titlesec}
\usepackage{tikz}
\usetikzlibrary{calc}
\usetikzlibrary{arrows}
\usetikzlibrary{decorations.markings}
\usepackage{geometry}

\newtheorem{theorem}{Theorem}[section]
\numberwithin{equation}{theorem}
\newtheorem{lemma}[theorem]{Lemma}

\newtheorem{corollary}[theorem]{Corollary}

\theoremstyle{definition}
\newtheorem{definition}[theorem]{Definition}
\newtheorem{example}[theorem]{Example}

\newtheorem{remark}[theorem]{Remark}

\theoremstyle{conjecture}

\newtheorem{question}[theorem]{Question}

\newcommand{\Ass}{\operatorname{Ass}}
\newcommand{\im}{\operatorname{im}}

\newcommand{\Spec}{\operatorname{Spec}}

\newcommand{\ann}{\operatorname{ann}}
\newcommand{\rank}{\operatorname{rank}}

\newcommand{\cd}{\operatorname{cd}}

\newcommand{\Ht}{\operatorname{ht}}

\newcommand{\V}{\operatorname{V}}

\newcommand{\Cone}{\operatorname{Cone}}

\newcommand{\Ext}{\operatorname{Ext}}
\newcommand{\Supp}{\operatorname{Supp}}

\newcommand{\Tor}{\operatorname{Tor}}
\newcommand{\Hom}{\operatorname{Hom}}

\newcommand{\depth}{\operatorname{depth}}

\newcommand{\coker}{\operatorname{coker}}

\newcommand{\Max}{\operatorname{Max}}

\newcommand{\lo}{\longrightarrow}

\newcommand{\fp}{\frak{p}}

\newcommand{\fa}{\frak{a}}
\newcommand{\fb}{\frak{b}}

\newcommand{\suchthat}{\;\ifnum\currentgrouptype=16 \middle\fi|\;}

\newenvironment{prf}[1][Proof]{\begin{proof}[\bf #1]}{\end{proof}}

\makeatletter
\newcommand{\holim@}[2]{%
  \vtop{\m@th\ialign{##\cr
    \hfil$#1\operator@font holim$\hfil\cr
    \noalign{\nointerlineskip\kern1.5\ex@}#2\cr
    \noalign{\nointerlineskip\kern-\ex@}\cr}}%
}
\newcommand{\holim}{%
  \mathop{\mathpalette\holim@{\rightarrowfill@\textstyle}}\nmlimits@
}
\makeatother

\makeatletter
\def\@secnumfont{\bfseries}
\def\section{\@startsection{section}{1}%
  \z@{.7\linespacing\@plus\linespacing}{.5\linespacing}%
  {\normalfont\Large\bfseries\filcenter}}
\def\subsection{\@startsection{subsection}{2}%
  \z@{.5\linespacing\@plus.7\linespacing}{-.5em}%
  {\normalfont\large\bfseries}}
\makeatother

\begin{document}

\author[K. Divaani-Aazar, H. Faridian and M. Tousi]{Kamran Divaani-Aazar, Hossein Faridian
and Massoud Tousi}

\title[Stable Under Specialization Sets and Cofiniteness]
{Stable Under Specialization Sets and Cofiniteness}

\address{K. Divaani-Aazar, Department of Mathematics, Alzahra University, Vanak, Post Code 19834, Tehran, Iran-and-School
of Mathematics, Institute for Research in Fundamental Sciences (IPM), P.O. Box 19395-5746, Tehran, Iran.}
\email{kdivaani@ipm.ir}

\address{H. Faridian, Department of Mathematics, Shahid Beheshti University, G.C., Evin, Tehran, Iran, Zip Code 1983963113.}
\email{h.faridian@yahoo.com}

\address{M. Tousi, Department of Mathematics, Shahid Beheshti University, G.C., Evin, Tehran, Iran, P.O. Box 19395-5746.}
\email{mtousi@ipm.ir}


\keywords {Abelian category; cofinite module; cohomological dimension; derived category; local cohomology module; stable
under specialization set.\\
The research of the first author is supported by a grant from IPM (No. 95130212).}

\begin{abstract}
Let $R$ be a commutative noetherian ring, and $\mathcal{Z}$ a stable under specialization subset of $\Spec(R)$.
We introduce a notion of $\mathcal{Z}$-cofiniteness and study its main properties. In the case $\dim(\mathcal{Z})\leq 1$,
or $\dim(R)\leq 2$, or $R$ is semilocal with $\cd(\mathcal{Z},R) \leq 1$,  we show that the category of $\mathcal{Z}$-cofinite
$R$-modules is abelian. Also, in each of these cases, we prove that the local cohomology module $H^{i}_{\mathcal{Z}}(X)$ is
$\mathcal{Z}$-cofinite for every homologically left-bounded $R$-complex $X$ whose homology modules are finitely generated
and every $i \in \mathbb{Z}$.
\end{abstract}

\maketitle

\tableofcontents

\section{Introduction}

\sloppy

Throughout this paper, $R$ denotes a commutative noetherian ring with identity and $\mathcal{M}(R)$ flags the category of
$R$-modules.

In his algebraic geometry seminars of 1961-2, Grothendieck founded the theory of local cohomology and raised, along the way,
a few questions on the finiteness properties of the local cohomology modules; see \cite[Conjectures 1.1 and 1.2]{Gr}. He specifically
asked whether the $R$-modules $\Hom_{R}\left(R/\mathfrak{a},H^{i}_{\mathfrak{a}}(M) \right)$ were finitely generated for every ideal
$\mathfrak{a}$ of $R$ and every finitely generated $R$-module $M$, which had been answered affirmatively in the same seminar when
$(R,\mathfrak{m})$ is local and $\mathfrak{a}=\mathfrak{m}$. In 1969, Hartshorne provided a counterexample in \cite[Section 3]{Ha1},
to show that this question does not have an affirmative answer in general.
For a given ideal $\mathfrak{a}$ of $R$, Hartshorne defined an $R$-module $M$ to be $\mathfrak{a}$-\textit{cofinite} if $\Supp_{R}(M)
\subseteq \V(\mathfrak{a})$ and $\Ext^{i}_{R}(R/\mathfrak{a},M)$ is finitely generated for every $i \geq 0$, and accordingly posed the
following questions:

\begin{question} \label{1.1}
Are the local cohomology modules $H^{i}_{\mathfrak{a}}(M)$, $\mathfrak{a}$-cofinite for every finitely generated $R$-module $M$
and every $i \geq 0$?
\end{question}

\begin{question} \label{1.2}
Is the category $\mathcal{M}(R,\mathfrak{a})_{cof}$ of $\mathfrak{a}$-cofinite $R$-modules an abelian subcategory of $\mathcal{M}(R)$?
\end{question}

By providing the following counterexample, he showed that the answers to these questions are negative in general.

\begin{example} \label{1.3}
Let $k$ be a field, $R=k[X,Y][[U,V]]$, $\mathfrak{a}=(U,V)$, $\mathfrak{p}=(XV+YU)$, and $T=R/\mathfrak{p}$. Then $R$ is a regular domain
of dimension $4$ and $T$ is a non-regular domain. It is shown that $\Hom_{R}\left(R/\mathfrak{a},H^{2}_{\mathfrak{a}}(T)\right)$ is not
finitely generated, so in particular, $H^{2}_{\mathfrak{a}}(T)$ is not $\mathfrak{a}$-cofinite. This takes care of Question \ref{1.1}.
Furthermore, there is an exact sequence $$0 \rightarrow H^{1}_{\mathfrak{a}}(T) \rightarrow H^{2}_{\mathfrak{a}}(R) \rightarrow H^{2}_{\mathfrak{a}}(R)
\rightarrow H^{2}_{\mathfrak{a}}(T)\rightarrow 0.$$ The local cohomology module $H^{2}_{\mathfrak{a}}(R)$ turns out to be $\mathfrak{a}$-cofinite,
whereas $H^{2}_{\mathfrak{a}}(T)$ is not $\mathfrak{a}$-cofinite, answering Question \ref{1.2}.
\end{example}

Hartshorne further established affirmative answers to these questions in the case where $\mathfrak{a}$ is a principal ideal generated by a
nonzerodivisor and $R$ is an $\mathfrak{a}$-adically complete regular ring of finite Krull dimension, and also in the case where $\mathfrak{a}$
is a prime ideal with $\dim(R/\mathfrak{a})=1$ and $R$ is a complete regular local ring; see \cite[Propositions 6.1
and 7.6, and Corollaries 6.3 and  7.7]{Ha1}.

In the following years, Hartshorne's results on Questions \ref{1.1} and \ref{1.2}, were systematically extended and polished by commutative algebra
practitioners in several stages to take the following full-fledged culminating form.

\begin{theorem} \label{1.4}
Let $\mathfrak{a}$ be an ideal of $R$, and $M$ a finitely generated $R$-module. Suppose that either of the following conditions are satisfied:
\begin{enumerate}
\item[(i)] $\cd(\mathfrak{a},R) \leq 1$, or
\item[(ii)] $\dim(R/\mathfrak{a})\leq 1$, or
\item[(iii)] $\dim(R) \leq 2$,
\end{enumerate}
Then $H^{i}_{\mathfrak{a}}(M)$ is $\mathfrak{a}$-cofinite for every $i \geq 0$, and $\mathcal{M}(R,\mathfrak{a})_{cof}$ is an abelian subcategory
of $\mathcal{M}(R)$.
\end{theorem}

For (i), see \cite[Corollary 3.14]{Me2} and \cite[Theorem 2.2 (ii)]{DFT}. For (ii), see \cite[Theorem 2.6 and Corollary 2.12]{Me1}, \cite[Corollary 2.8]{BNS},
and \cite[Corollary 2.7]{BN}. Finally for (iii), see \cite[Theorems 7.10 and 7.4]{Me2}.

For an $\frak{a}$-cofinite $R$-module $M$, it is known that the set $\Ass_R(M)$ is finite, and all its Bass numbers and Betti numbers with respect to every prime ideal of $R$ are also finite. Probing such finiteness properties has been a high-profile problem in commutative algebra; see e.g. \cite{HK}, \cite{HS} and \cite{Ly}.

Several authors have strived to extend the results of Theorem \ref{1.4} to generalized local cohomology modules. However, it is folklore that all the generalizations
$H_{\varphi}^{i}(M)$, $H_{\mathcal{Z}}^{i}(M)$, $H_{\mathcal{\mathfrak{a},\mathfrak{b}}}^{i}(M)$, $H_{\mathcal{\mathfrak{a}}}^{i}(M,N)$, $H_{\varphi}^{i}(M,N)$, and $H_{\mathcal{\mathfrak{a},\mathfrak{b}}}^{i}(M,N)$ of the local cohomology module $H_{\mathcal{\mathfrak{a}}}^{i}(M)$ of an $R$-module $M$, are special cases of the
local cohomology module $H_{\mathcal{Z}}^{i}(X)$ of an $R$-complex $X$ with support in a stable under specialization subset $\mathcal{Z}$ of $\Spec(R)$; see Remark
\ref{4.9}. Therefore, any established result on $H_{\mathcal{Z}}^{i}(X)$ encompasses all the previously known results on each of these local cohomology modules.
In this direction, we aspire to define the general notion of $\mathcal{Z}$-cofiniteness and extend Theorem \ref{1.4} to $H^{i}_{\mathcal{Z}}(X)$. We specifically
obtain the following results; see Theorems \ref{3.11}, and \ref{4.8}.

\begin{theorem} \label{1.5}
Let $\mathcal{Z}$ be a stable under specialization subset of $\Spec (R)$ such that either $R$ is semilocal with $\cd(\mathcal{Z},R) \leq 1$, or $\dim(\mathcal{Z})
\leq 1$, or $\dim(R) \leq 2$. Then $\mathcal{M}(R,\mathcal{Z})_{cof}$ is an abelian subcategory of $\mathcal{M}(R)$.
\end{theorem}

\begin{theorem} \label{1.6}
Let $\mathcal{Z}$ a stable under specialization subset of $\Spec(R)$, and $X$ a homologically left-bounded $R$-complex $X$ with finitely generated homology modules.
Assume that either $R$ is semilocal with $\cd(\mathcal{Z},R)\leq 1$, or $\dim(\mathcal{Z})\leq 1$, or $\dim\left(\Supp_{R}(X)\right) \leq 2$. Then $H^{i}_{\mathcal{Z}}(X)$ is $\mathcal{Z}$-cofinite for every $i \in \mathbb{Z}$.
\end{theorem}

\section{Preliminaries}

In this section, we first present some background material on complexes which will be used in the rest of the work. For more information, refer to \cite{AF}, \cite{CFH}, \cite{Ha2}, \cite{Li}, and \cite{Sp}. In what follows, $\mathcal{C}(R)$ denotes the category of $R$-complexes.

The derived category $\mathcal{D}(R)$ is defined as the localization of the homotopy category $\mathcal{K}(R)$ with respect to the multiplicative system of quasi-isomorphisms. Simply put, an object in $\mathcal{D}(R)$ is an $R$-complex $X$ displayed in the standard homological style
$$X= \cdots \rightarrow X_{i+1} \xrightarrow {\partial^{X}_{i+1}} X_{i} \xrightarrow {\partial^{X}_{i}} X_{i-1} \rightarrow \cdots,$$
and a morphism $\varphi:X\rightarrow Y$ in $\mathcal{D}(R)$ is given by the equivalence class of a pair $(f,g)$ of morphisms
$X \xleftarrow{g} U \xrightarrow{f} Y$ in $\mathcal{C}(R)$ with $g$ a quasi-isomorphism, under the equivalence relation that identifies two such pairs $(f,g)$ and $(f^{\prime},g^{\prime})$, whenever there is a diagram in $\mathcal{C}(R)$ as follows which commutes up to homotopy:
\[\begin{tikzpicture}[every node/.style={midway}]
  \matrix[column sep={3em}, row sep={3em}]
  {\node(1) {$$}; & \node(2) {$U$}; & \node(3) {$$};\\
  \node(4) {$X$}; & \node(5) {$V$}; & \node(6) {$Y$};\\
  \node(7) {$$}; & \node(8) {$U^{\prime}$}; & \node(9) {$$};\\};
  \draw[decoration={markings,mark=at position 1 with {\arrow[scale=1.5]{>}}},postaction={decorate},shorten >=0.5pt] (2) -- (4) node[anchor=south] {$g$} node[anchor=west] {$\simeq$};
  \draw[decoration={markings,mark=at position 1 with {\arrow[scale=1.5]{>}}},postaction={decorate},shorten >=0.5pt] (2) -- (6) node[anchor=south] {$f$};
  \draw[decoration={markings,mark=at position 1 with {\arrow[scale=1.5]{>}}},postaction={decorate},shorten >=0.5pt] (8) -- (4) node[anchor=north] {$g^{\prime}$} node[anchor=west] {$\simeq$};
  \draw[decoration={markings,mark=at position 1 with {\arrow[scale=1.5]{>}}},postaction={decorate},shorten >=0.5pt] (8) -- (6) node[anchor=north] {$f^{\prime}$};
  \draw[decoration={markings,mark=at position 1 with {\arrow[scale=1.5]{>}}},postaction={decorate},shorten >=0.5pt] (5) -- (2) node[anchor=east] {$$};
  \draw[decoration={markings,mark=at position 1 with {\arrow[scale=1.5]{>}}},postaction={decorate},shorten >=0.5pt] (5) -- (4) node[anchor=south] {$\simeq$};
  \draw[decoration={markings,mark=at position 1 with {\arrow[scale=1.5]{>}}},postaction={decorate},shorten >=0.5pt] (5) -- (6) node[anchor=east] {$$};
  \draw[decoration={markings,mark=at position 1 with {\arrow[scale=1.5]{>}}},postaction={decorate},shorten >=0.5pt] (5) -- (8) node[anchor=east] {$$};
\end{tikzpicture}\]
The isomorphisms in $\mathcal{D}(R)$ are marked by the symbol $\simeq$.

The derived category $\mathcal{D}(R)$ is triangulated. A distinguished triangle in $\mathcal{D}(R)$ is a triangle that is isomorphic to a triangle of the form
$$X \xrightarrow {\mathfrak{L}(f)} Y \xrightarrow{\mathfrak{L}(\varepsilon)} \Cone(f) \xrightarrow{\mathfrak{L}(\varpi)} \Sigma X,$$
for some morphism $f:X \rightarrow Y$ in $\mathcal{C}(R)$ with the mapping cone sequence
$$0 \rightarrow Y \xrightarrow{\varepsilon} \Cone(f) \xrightarrow{\varpi} \Sigma X \rightarrow 0,$$
in which $\mathfrak{L}:\mathcal{C}(R) \rightarrow \mathcal{D}(R)$ is the canonical functor that is defined as $\mathfrak{L}(X)=X$ for every $R$-complex $X$, and $\mathfrak{L}(f)=\varphi$ where $\varphi$ is represented by the morphisms $X \xleftarrow{1^{X}} X \xrightarrow{f} Y$ in $\mathcal{C}(R)$.

We let $\mathcal{D}_{\sqsubset}(R)$ (res. $\mathcal{D}_{\sqsupset}(R)$) denote the full subcategory of $\mathcal{D}(R)$ consisting of $R$-complexes $X$ with $H_{i}(X)=0$ for $i \gg 0$ (res. $i \ll 0$), and  $D_{\square}(R):=\mathcal{D}_{\sqsubset}(R)\cap \mathcal{D}_{\sqsupset}(R)$. We further let $\mathcal{D}^{f}(R)$ denote the full subcategory of $\mathcal{D}(R)$ consisting of $R$-complexes $X$ with finitely generated homology modules. We also feel free to use any combination of the subscripts and the superscript as in $\mathcal{D}^{f}_{\square}(R)$, with the obvious meaning of the intersection of the two subcategories involved. Given an $R$-complex $X$, the standard notions
$$\sup(X) = \sup \left\{i \in \mathbb{Z} \suchthat H_{i}(X) \neq 0 \right\}$$
and
$$\inf(X) = \inf \left\{i \in \mathbb{Z} \suchthat H_{i}(X) \neq 0 \right\}$$
are frequently used, with the convention that $\sup(\emptyset) =- \infty$ and $\inf(\emptyset) = +\infty$.

An $R$-complex $P$ of projective modules is said to be semi-projective if the functor $\Hom_{R}(P,-)$ preserves quasi-isomorphisms. By a semi-projective resolution of an $R$-complex $X$, we mean a quasi-isomorphism $P\xrightarrow {\simeq} X$ in which $P$ is a semi-projective $R$-complex. Dually, an $R$-complex
$I$ of injective modules is said to be semi-injective if the functor $\Hom_{R}(-,I)$ preserves quasi-isomorphisms. By a semi-injective resolution of an $R$-complex $X$, we mean a quasi-isomorphism $X\xrightarrow {\simeq} I$ in which $I$ is a semi-injective $R$-complex. Semi-projective and semi-injective resolutions exist for any $R$-complex; see \cite[Theorems 5.2.13 and 5.3.18]{CFH}. Moreover, any right-bounded $R$-complex of projective modules is semi-projective, and any left-bounded $R$-complex of injective modules is semi-injective; see \cite[ Examples 5.2.7 and 5.3.11]{CFH}.

Let $X$ and $Y$ be two $R$-complexes. Then each of the functors $\Hom_{R}(X,-)$ and $\Hom_{R}(-,Y)$ on $\mathcal{C}(R)$ preserves homotopy equivalences, and thus enjoys a right total derived functor on $\mathcal{D}(R)$ \cite[Theorem 7.1.14]{CFH}, together with a balance property, in the sense that ${\bf R}\Hom_{R}(X,Y)$ can be computed by
$${\bf R}\Hom_{R}(X,Y)\simeq \Hom_{R}(P,Y) \simeq \Hom_{R}(X,I),$$
where $P\xrightarrow {\simeq} X$ is any semi-projective resolution of $X$, and $Y\xrightarrow {\simeq} I$ is any semi-injective resolution of $Y$.   Moreover, we let $\Ext^{i}_{R}(X,Y):=H_{-i}\left({\bf R}\Hom_{R}(X,Y)\right)$ for every $i \in \mathbb{Z}$.

Likewise, each of the functors $X\otimes_{R}-$ and $-\otimes_{R}Y$ on $\mathcal{C}(R)$ preserves homotopy equivalences, and thus enjoys a left total derived functor on $\mathcal{D}(R)$ \cite[Theorem 7.2.20]{CFH}, together with a balance property, in the sense that $X\otimes_{R}^{\bf L}Y$ can be computed by
$$X\otimes_{R}^{\bf L}Y \simeq P\otimes_{R}Y \simeq X\otimes_{R}Q,$$
where $P\xrightarrow {\simeq} X$ is any semi-projective resolution of $X$, and $Q\xrightarrow {\simeq} Y$ is any semi-projective resolution of $Y$. Moreover, we let $\Tor^{R}_{i}(X,Y):=H_{i}\left(X\otimes_{R}^{\bf L}Y\right)$
for every $i \in \mathbb{Z}$.

We next turn to the notion of a stable under specialization set. A subset $\mathcal{Z}$ of $\Spec(R)$ is said to be stable under specialization if $V(\fp)\subseteq \mathcal{Z}$ for every $\fp\in \mathcal{Z}$. For such a subset
$\mathcal{Z}$, we set
$F(\mathcal{Z}):=\left\{\mathfrak{a}\lhd R \suchthat V(\mathfrak{a})\subseteq \mathcal{Z}\right\}$.
If $M$ is an $R$-module, then clearly $\Supp_{R}(M)$ is a stable under specialization subset of $\Spec(R)$. Conversely, given any stable under specialization subset $\mathcal{Z}$ of $\Spec(R)$, one readily checks out that $\mathcal{Z} = \Supp_{R}\left(\bigoplus_{\mathfrak{a} \in F(\mathcal{Z})} R/\mathfrak{a}\right)$.
In particular, $\V(\mathfrak{a})$ for an ideal $\mathfrak{a}$ of $R$, and any subset of $\Max(R)$  are stable under specialization subsets of $\Spec(R)$.

We finally recall the definition of the most general local cohomology functor. Given a stable under specialization subset $\mathcal{Z}$ of $\Spec(R)$, we let
$$\Gamma_{\mathcal{Z}}(M) := \left\{x\in M \suchthat \Supp_{R}(Rx)\subseteq \mathcal{Z}\right\}$$
for an $R$-module $M$, and $\Gamma_{\mathcal{Z}}(f) := f|_{\Gamma_{\mathcal{Z}}(M)}$ for an $R$-homomorphism $f:M\rightarrow N$. This provides us with the so-called $\mathcal{Z}$-torsion functor $\Gamma_{\mathcal{Z}}(-)$ on $\mathcal{M}(R)$, which extends, by termwise action, to a functor on $\mathcal{C}(R)$. The extended functor clearly preserves homotopy equivalences. Therefore, it enjoys a right total derived functor ${\bf R}\Gamma_{\mathcal{Z}}(-)$ on $\mathcal{D}(R)$ \cite[Definition 6.6.12]{CFH}, that can be computed by ${\bf R}\Gamma_{\mathcal{Z}}(X)\simeq \Gamma_{\mathcal{Z}}(I)$, where $X \xrightarrow {\simeq} I$ is any semi-injective resolution of $X$. Besides, we define the $i$th local cohomology module of $X$ with support in $\mathcal{Z}$ as $H^{i}_{\mathcal{Z}}(X):= H_{-i}\left({\bf R}\Gamma_{\mathcal{Z}}(X)\right)$ for every $i\in\mathbb{Z}$. It is obvious that upon setting $\mathcal{Z}=\V(\mathfrak{a})$ for some ideal $\mathfrak{a}$ of $R$, we recover the usual local cohomology module with respect to $\mathfrak{a}$.

It is straightforward to see that the set $F(\mathcal{Z})$ is a directed poset under reverse inclusion. Let $X$ be an $R$-complex and $X \xrightarrow {\simeq} I$ a semi-injective resolution of $X$. Then one can see by inspection that
$$\Gamma_{\mathcal{Z}}(I_{i})=\bigcup_{\mathfrak{a}\in F(\mathcal{Z})}\Gamma_{\mathfrak{a}}(I_{i}) \cong \underset{\fa\in F(\mathcal{Z})}{\varinjlim}\Gamma_{\mathfrak{a}}(I_{i})$$
for every $i\in \mathbb{Z}$, which in turn implies that $\Gamma_{\mathcal{Z}}(I) \cong \underset{\fa\in F(\mathcal{Z})}{\varinjlim}\Gamma_{\mathfrak{a}}(I)$. Therefore, we have

\begin{equation*}
\begin{split}
H_{\mathcal{Z}}^i(X) & \cong H_{-i}\left({\bf R}\Gamma_{\mathcal{Z}}(X)\right) \\
 & \cong H_{-i}\left(\Gamma_{\mathcal{Z}}(I)\right) \\
 & \cong H_{-i}\left(\underset{\fa\in F(\mathcal{Z})}{\varinjlim}\Gamma_{\mathfrak{a}}(I)\right) \\
 & \cong \underset{\fa\in F(\mathcal{Z})}{\varinjlim} H_{-i}\left(\Gamma_{\mathfrak{a}}(I)\right) \\
 & \cong \underset{\fa\in F(\mathcal{Z})}{\varinjlim} H_{-i}\left({\bf R}\Gamma_{\mathfrak{a}}(X)\right) \\
 & \cong \underset{\fa\in F(\mathcal{Z})}{\varinjlim}H_{\mathfrak{a}}^i(X)
\end{split}
\end{equation*}
for every $i \in \mathbb{Z}$.

For a stable under specialization subset $\mathcal{Z}$ of $\Spec(R)$, we define the dimension of $\mathcal{Z}$ as $$\dim(\mathcal{Z}):=\sup \left\{\dim(R/\mathfrak{a}) \suchthat \mathfrak{a}\in F(\mathcal{Z})\right\}.$$ Also, we define the cohomological dimension of an $R$-complex $X$ with respect to $\mathcal{Z}$ as $$\cd(\mathcal{Z},X):=\sup \left\{i
\in \mathbb{Z} \suchthat H^{i}_{\mathcal{Z}}(X)\neq 0 \right\}.$$

Now, we are ready to define the general notion of $\mathcal{Z}$-cofiniteness. Recall that the support of an $R$-complex $X$ is defined to be $\Supp_{R}(X)=\bigcup_{i \in \mathbb{Z}}\Supp_{R}\left(H_{i}(X)\right)$.

\begin{definition} \label{2.1}
Let $\mathcal{Z}$ be a stable under specialization subset of $\Spec(R)$. An $R$-complex $X\in \mathcal{D}(R)$ is said to be $\mathcal{Z}$\textit{-cofinite} if $\Supp_{R}(X)\subseteq \mathcal{Z}$ and ${\bf R}\Hom_{R}(R/\mathfrak{a},X)\in \mathcal{D}^{f}(R)$ for every $\mathfrak{a}\in F(\mathcal{Z})$.
\end{definition}

We denote the full subcategory of $\mathcal{M}(R)$ consisting of $\mathcal{Z}$-cofinite $R$-modules by $\mathcal{M}(R,\mathcal{Z})_{cof}$. The next result lays on some characterizations of $\mathcal{Z}$-cofinite complexes.

\begin{lemma} \label{2.2}
Let $\mathcal{Z}$ be a stable under specialization subset of $\Spec(R)$, and $X \in \mathcal{D}(R)$ with $\Supp_{R}(X)\subseteq \mathcal{Z}$. Consider the following conditions:
\begin{enumerate}
\item[(a)] $X$ is $\mathcal{Z}$-cofinite.
\item[(b)] ${\bf R}\Hom_{R}(Y,X)\in \mathcal{D}^{f}(R)$ for every $Y \in \mathcal{D}^{f}_{\square}(R)$ with $\Supp_{R}(Y)\subseteq \mathcal{Z}$.
\item[(c)] $(R/\mathfrak{a})\otimes_{R}^{\bf L} X\in \mathcal{D}^{f}(R)$ for every $\mathfrak{a}\in F(\mathcal{Z})$.
\item[(d)] $Y\otimes_{R}^{\bf L} X\in \mathcal{D}^{f}(R)$ for every $Y \in \mathcal{D}^{f}_{\square}(R)$ with $\Supp_{R}(Y)\subseteq \mathcal{Z}$.
\end{enumerate}
Then the following assertions hold:
\begin{enumerate}
\item[(i)] If $X \in \mathcal{D}_{\sqsubset}(R)$, then (a) and (b) are equivalent.
\item[(ii)] If $X \in \mathcal{D}_{\sqsupset}(R)$, then (c) and (d) are equivalent.
\item[(iii)] If $X \in \mathcal{D}_{\square}(R)$, then all the assertions are equivalent.
\end{enumerate}
\end{lemma}

\begin{prf}
(i): Suppose that (a) holds and $Y \in \mathcal{D}^{f}_{\square}(R)$ with $\Supp_{R}(Y)\subseteq \mathcal{Z}$. Let $\mathfrak{a}=\ann_{R}\left(\oplus_{i\in \mathbb{Z}}H_{i}(Y)\right)$.
Then $V(\mathfrak{a})= \Supp_{R}(Y) \subseteq \mathcal{Z}$, so $\mathfrak{a}\in F(\mathcal{Z})$. Therefore, ${\bf R}\Hom_{R}(R/\mathfrak{a},X)\in\mathcal{D}^{f}(R)$.
Now, \cite[Proposition 7.2]{WW} implies that ${\bf R}\Hom_{R}(Y,X)\in \mathcal{D}^{f}(R)$. The converse is clear.

(ii): Similar to (i) using \cite[Proposition 7.1]{WW}.

(iii): Fix $\mathfrak{a}\in F(\mathcal{Z})$. Then \cite[Proposition 7.4]{WW} yields that ${\bf R}\Hom_{R}(R/\mathfrak{a},X)\in \mathcal{D}^{f}(R)$ if and only if $(R/\mathfrak{a})\otimes_{R}^{\bf L} X\in \mathcal{D}^{f}(R)$.
\end{prf}

We collect some basic properties of $\mathcal{Z}$-cofinite $R$-complexes in the following result. Its first part indicates that in the case where $\mathcal{Z}= V(\mathfrak{a})$, our definition of $\mathcal{Z}$-cofiniteness coincides with the usual notion of $\mathfrak{a}$-cofiniteness.

\begin{lemma} \label{2.3}
Let $\mathcal{Z}$ be a stable under specialization subset of $\Spec(R)$ and $X\in \mathcal{D}_{\sqsubset}(R)$. Then the following assertions hold:
\begin{enumerate}
\item[(i)] If $\mathfrak{a}$ is an ideal of $R$, then $X$ is $V(\mathfrak{a})$-cofinite if and only if $\Supp_{R}(X)\subseteq V(\mathfrak{a})$ and ${\bf R}\Hom_{R}(R/\mathfrak{a},X)\in \mathcal{D}^{f}(R)$.
\item[(ii)] If $H_{i}(X)$ is $\mathcal{Z}$-cofinite for every $i\in \mathbb{Z}$, then $X$ is $\mathcal{Z}$-cofinite.
\item[(iii)] If $X\in \mathcal{D}^{f}_{\sqsubset}(R)$, then ${\bf R}\Gamma_{\mathcal{Z}}(X)$ is $\mathcal{Z}$-cofinite.
\item[(iv)] If $X$ is $\mathcal{Z}$-cofinite, then the Bass number $\mu^{i}(\mathfrak{p},X)$ is finite for every $\mathfrak{p}\in \Spec (R)$ and every $i \in \mathbb{Z}$.
\item[(v)] If $X\in \mathcal{D}_{\square}(R)$ and $X$ is $\mathcal{Z}$-cofinite, then the Betti number $\beta_{i}(\mathfrak{p},X)$ is finite for every $\mathfrak{p}\in \Spec (R)$ and every $i \in \mathbb{Z}$.
\end{enumerate}
\end{lemma}

\begin{prf}
(i): Suppose that $\Supp_{R}(X)\subseteq V(\mathfrak{a})$ and ${\bf R}\Hom_{R}(R/\mathfrak{a},X)\in \mathcal{D}^{f}(R)$. Let $\mathfrak{b}\in F\left(V(\mathfrak{a})\right)$. It follows that $\Supp_{R}(R/\mathfrak{b}) \subseteq V(\mathfrak{a})$. Now, by \cite[Proposition 7.2]{WW} we are through. The converse is clear.

(ii): Since $H_{i}(X)$ is $\mathcal{Z}$-cofinite for every $i\in \mathbb{Z}$, we have
$$\Supp_{R}(X)=\bigcup_{i\in \mathbb{Z}}\Supp_{R}\left(H_{i}(X)\right) \subseteq \mathcal{Z}.$$
Let $\mathfrak{a}\in F(\mathcal{Z})$.
The spectral sequence
$$E_{p,q}^{2}=\Ext_{R}^{p}\left(R/\mathfrak{a},H_{-q}(X)\right) \underset{p}\Rightarrow \Ext_{R}^{p+q}(R/\mathfrak{a},X)$$
from the proof of \cite[Proposition 6.2]{Ha1}, together with the assumption that $E_{p,q}^{2}$ is finitely generated for every $p,q\in \mathbb{Z}$, conspire to imply that $\Ext_{R}^{p+q}(R/\mathfrak{a},X)$ is finitely generated, i.e. $X$ is $\mathcal{Z}$-cofinite.

(iii): It is clear that $\Supp_{R}\left({\bf R}\Gamma_{\mathcal{Z}}(X)\right)\subseteq \mathcal{Z}$. Since $X$ is homologically left-bounded, there exists a left-bounded semi-injective resolution $X \xrightarrow {\simeq} I$ of $X$ invoking \cite[Theorem 5.3.26]{CFH}.
As $I$ is an $R$-complex of injective modules, and $\Gamma_{\mathcal{Z}}(I_{i}) \cong \underset{\fa\in F(\mathcal{Z})}{\varinjlim}\Gamma_{\mathfrak{a}}(I_{i})$
for every $i \in \mathbb{Z}$, we see that $\Gamma_{\mathcal{Z}}(I)$ is a left-bounded $R$-complex of injective modules, and thus $\Gamma_{\mathcal{Z}}(I)$ is semi-injective.

Fix $\mathfrak{a}\in F(\mathcal{Z})$. For every $R$-module $\Hom_{R}\left(R/\mathfrak{a},\Gamma_{\mathcal{Z}}(I_{i})\right)$ in the $R$-complex $\Hom_{R}\left(R/\mathfrak{a},\Gamma_{\mathcal{Z}}(I)\right)$, we have
$$\Hom_{R}\left(R/\mathfrak{a},\Gamma_{\mathcal{Z}}(I_{i})\right) \cong \Hom_{R}\left(R/\mathfrak{a},\Gamma_{\mathfrak{a}}(I_{i})\right) \cong \Hom_{R}\left(R/\mathfrak{a},I_{i}\right).$$
Indeed, if $N$ is a submodule of $M$ containing $\Gamma_{\mathfrak{a}}(M)$, then every $R$-homomorphism $f: R/\mathfrak{a} \rightarrow N$ has its image in $\Gamma_{\mathfrak{a}}(M)$.
Hence one has
\begin{equation*}
\begin{split}
{\bf R}\Hom_{R}\left(R/\mathfrak{a},{\bf R}\Gamma_{\mathcal{Z}}(X)\right) & \simeq {\bf R}\Hom_{R}\left(R/\mathfrak{a},\Gamma_{\mathcal{Z}}(I)\right) \\
 & \simeq \Hom_{R}\left(R/\mathfrak{a},\Gamma_{\mathcal{Z}}(I)\right) \\
 & \simeq \Hom_{R}(R/\mathfrak{a},I) \\
 & \simeq {\bf R}\Hom_{R}(R/\mathfrak{a},X).
\end{split}
\end{equation*}
But $X\in \mathcal{D}^{f}_{\sqsubset}(R)$, so ${\bf R}\Hom_{R}(R/\mathfrak{a},X) \in \mathcal{D}^{f}_{\sqsubset}(R)$, and thus the assertion follows.

(iv): For every $i\in \mathbb{Z}$, we have
$$\mu^{i}(\mathfrak{p},X):= \rank_{R_{\mathfrak{p}}/\mathfrak{p}R_{\mathfrak{p}}}\left(\Ext_{R_{\mathfrak{p}}}^{i}\left(R_{\mathfrak{p}}/\mathfrak{p}R_{\mathfrak{p}},
X_{\mathfrak{p}}\right)\right).$$
If $\mathfrak{p}\not\in \mathcal{Z}$,
then $\mathfrak{p}\not\in \Supp_{R}(X)$, so $\mu^{i}(\mathfrak{p},X)=0$. If $\mathfrak{p}\in \mathcal{Z}$, then $\V(\mathfrak{p}) \subseteq \mathcal{Z}$, so by definition, ${\bf R}\Hom_{R}(R/\mathfrak{p},X)\in \mathcal{D}^{f}(R)$, whence
$\mu^{i}(\mathfrak{p},X)< \infty$.

(v): For every $i\in \mathbb{Z}$, we have
$$\beta_{i}(\mathfrak{p},X):= \rank_{R_{\mathfrak{p}}/\mathfrak{p}R_{\mathfrak{p}}}\left(\Tor_{i}^{R_{\mathfrak{p}}}\left(R_{\mathfrak{p}}/\mathfrak{p}R_{\mathfrak{p}},
X_{\mathfrak{p}}\right)\right).$$
If $\mathfrak{p}\not\in \mathcal{Z}$, then as in (iv), $\beta_{i}(\mathfrak{p},X)=0$. If $\mathfrak{p}\in
\mathcal{Z}$, then Lemma \ref{2.2} (iii) implies that $(R/\mathfrak{p})\otimes_{R}^{\bf L} X\in \mathcal{D}^{f}(R)$, thereby
$\beta_{i}(\mathfrak{p},X)<\infty$.
\end{prf}

\section{Proof of Theorem 1.5}

Our aim in this section is to prove Theorem \ref{1.5}; see Theorem \ref{3.11}. Corollary \ref{3.5} is the main ingredient in the proof of the first part of Theorem \ref{3.11}.
To prove it, we need Lemmas \ref{3.1}, \ref{3.2} and  \ref{3.4}.

\begin{lemma} \label{3.1}
Let $\mathcal{Z}$ be a stable under specialization subset of $\Spec(R)$. Let $M$ be an $R$-module such that $\Ass_{R}(M) \cap \mathcal{Z} \cap \Max(R)$ is a
finite set. Then the following assertions are equivalent:
\begin{enumerate}
\item[(i)] $(0:_{M}\mathfrak{a})$ is an artinian $R$-module for every $\mathfrak{a} \in F(\mathcal{Z})$.
\item[(ii)] $\Gamma_{\mathfrak{a}}(M)$ is an artinian $R$-module for every $\mathfrak{a} \in F(\mathcal{Z})$.
\item[(iii)] $\Gamma_{\mathcal{Z}}(M)$ is an artinian $R$-module.
\end{enumerate}
\end{lemma}

\begin{prf}
(i) $\Rightarrow$ (ii): Let $\mathfrak{a}\in F(\mathcal{Z})$. Then $\Gamma_{\mathfrak{a}}(M)$ is $\mathfrak{a}$-torsion, and
$(0:_{\Gamma_{\mathfrak{a}}(M)}\mathfrak{a}) = (0:_{M}\mathfrak{a})$
is artinian, so using Melkersson's Criterion \cite[Theorem 1.3]{Me3}, we conclude that $\Gamma_{\mathfrak{a}}(M)$ is artinian.

(ii) $\Rightarrow$ (iii): Set
$$\mathfrak{J}:= \bigcap_{\mathfrak{m}\in \Ass_{R}(M)\cap \mathcal{Z} \cap \Max(R)} \mathfrak{m}.$$
Since $\Ass_{R}(M)\cap \mathcal{Z} \cap \Max(R)$ is a finite set, we see that $V(\mathfrak{J})= \Ass_{R}(M)\cap \mathcal{Z} \cap \Max(R)$,
so that $\mathfrak{J} \in F(\mathcal{Z})$. Therefore, $\Gamma_{\mathfrak{J}}(M)$ is an artinian $R$-module.

Now, let $\mathfrak{a} \in F(\mathcal{Z})$. Then $\Gamma_{\mathfrak{a}}(M)$ is an artinian $R$-module by the assumption. Let $x \in \Gamma_{\mathfrak{a}}(M)$.
As $Rx$ is artinian, it turns out that
\begin{equation*}
\begin{split}
\Ass_{R}(Rx) & \subseteq \Ass_{R}(M) \cap V(\mathfrak{a}) \cap \Max(R) \\
 & \subseteq \Ass_{R}(M) \cap \mathcal{Z} \cap \Max(R) \\
 & = V(\mathfrak{J}).
\end{split}
\end{equation*}
Hence $x \in \Gamma_{\mathfrak{J}}(M)$. This yields that
$$\Gamma_{\mathcal{Z}}(M) = \bigcup_{\mathfrak{a}\in F(\mathcal{Z})} \Gamma_{\mathfrak{a}}(M) \subseteq \Gamma_{\mathfrak{J}}(M) \subseteq \Gamma_{\mathcal{Z}}(M),$$
thereby $\Gamma_{\mathcal{Z}}(M)=\Gamma_{\mathfrak{J}}(M)$ is artinian.

(iii) $\Rightarrow$ (i): Clear, since $(0:_{M}\mathfrak{a}) \subseteq \Gamma_{\mathcal{Z}}(M)$ for every $\mathfrak{a} \in F(\mathcal{Z})$.
\end{prf}

\begin{lemma} \label{3.2}
Let $\mathcal{Z}$ be a stable under specialization subset of $\Spec(R)$. Let $M$ be an $R$-module, and $r \geq 0$ an integer. Consider the following conditions:
\begin{enumerate}
\item[(a)] $H^{i}_{\mathcal{Z}}(M)$ is an artinian $R$-module for every $0 \leq i \leq r$.
\item[(b)] $\Ext^{i}_{R}(R/\mathfrak{a},M)$ is an artinian $R$-module for every $\mathfrak{a} \in F(\mathcal{Z})$ and for every $0 \leq i \leq r$.
\end{enumerate}
Then the following assertions hold:
\begin{enumerate}
\item[(i)] (a) implies (b).
\item[(ii)] If $\Supp_{R}(M)\cap \mathcal{Z} \cap \Max(R)$ is a finite set, then (a) and (b) are equivalent.
\end{enumerate}
\end{lemma}

\begin{prf}Let
$$I: 0\rightarrow I_{0} \xrightarrow {\partial^{I}_{0}} I_{-1} \xrightarrow {\partial^{I}_{-1}} I_{-2} \rightarrow \cdots$$
be a minimal injective resolution of $M$. Given any $\mathfrak{a}\in F(\mathcal{Z})$, consider the two $R$-complexes
$$\Hom_{R}(R/\mathfrak{a},I): 0\rightarrow \Hom_{R}(R/\mathfrak{a},I_{0}) \xrightarrow {\Hom_{R}\left(R/\mathfrak{a},\partial^{I}_{0}\right)} \Hom_{R}(R/\mathfrak{a},I_{-1}) $$$$ \xrightarrow {\Hom_{R}\left(R/\mathfrak{a},\partial^{I}_{-1}\right)} \Hom_{R}(R/\mathfrak{a},I_{-2}) \rightarrow \cdots,$$
and
$$\Gamma_{\mathcal{Z}}(I): 0\rightarrow \Gamma_{\mathcal{Z}}(I_{0}) \xrightarrow {\Gamma_{\mathcal{Z}}\left(\partial^{I}_{0}\right)} \Gamma_{\mathcal{Z}}(I_{-1}) \xrightarrow {\Gamma_{\mathcal{Z}}\left(\partial^{I}_{-1}\right)} \Gamma_{\mathcal{Z}}(I_{-2}) \rightarrow \cdots.$$
One can easily check that $\ker \left(\Hom_{R}\left(R/\mathfrak{a},\partial^{I}_{-i}\right)\right)$ is an essential submodule of $\Hom(R/\mathfrak{a},I_{-i})$, and $\ker \left(\Gamma_{\mathcal{Z}}\left(\partial^{I}_{-i}\right)\right)$ is an essential submodule of $\Gamma_{\mathcal{Z}}(I_{-i})$ for every $i \geq 0$.

Let $$X: 0 \rightarrow X_{0} \xrightarrow {\partial^{X}_{0}} X_{-1} \xrightarrow {\partial^{X}_{-1}} X_{-2} \rightarrow \cdots$$
be an $R$-complex such that $\ker \partial^{X}_{-i}$ is an essential submodule of $X_{-i}$ for every $i \geq 0$. For any given
$r\geq 0$, \cite[Lemma 5.4]{Me2} yields that $X_{-i}$ is an artinian $R$-module for every $0 \leq i \leq r$ if and only if $H_{-i}(X)$ is an artinian $R$-module for every $0 \leq i \leq r$. In the remainder of the proof, we apply this twice.

Now, we prove the following:

(i): Let $\mathfrak{a} \in F(\mathcal{Z})$. Applying the discussion above to the $R$-complex $\Gamma_{\mathcal{Z}}(I)$, we see that $\Gamma_{\mathcal{Z}}(I_{-i})$ is artinian for every $0 \leq i \leq r$. Since
$$\Hom_{R}(R/\mathfrak{a},I_{-i}) \cong (0:_{I_{-i}}\mathfrak{a}) \subseteq \Gamma_{\mathcal{Z}}(I_{-i}),$$
it is obvious that $\Hom_{R}(R/\mathfrak{a},I_{-i})$ is artinian for every $0 \leq i \leq r$.
This shows that $\Ext^{i}_{R}(R/\mathfrak{a},M)$ is artinian for every $0 \leq i \leq r$.

(ii): Another application of the discussion above to the $R$-complex $\Hom_{R}(R/\mathfrak{a},I)$, yields that $\Hom_{R}(R/\mathfrak{a},I_{-i})$ is artinian for every $\mathfrak{a} \in F(\mathcal{Z})$ and $0 \leq i \leq r$. Since $I$ is a minimal injective resolution of $M$, one can see that $\Supp_{R}(I_{-i}) \subseteq \Supp_{R}(M)$ for every $i \geq 0$. Therefore, we can use Lemma \ref{3.1} to deduce that $\Gamma_{\mathcal{Z}}(I_{-i})$ is artinian for every $0 \leq i \leq r$. A fortiori, $H^{i}_{\mathcal{Z}}(M)$ is artinian for
every $0 \leq i \leq r$.
\end{prf}

Part (i) of the next example shows that the finiteness condition on the sets $\Ass_{R}(M)\cap \mathcal{Z} \cap \Max(R)$ and $\Supp_{R}(M)\cap \mathcal{Z}\cap
\Max(R)$ cannot be removed from Lemmas \ref{3.1} and \ref{3.2}; respectively. Part (ii) of this example also demonstrates that unlike $\mathfrak{a}$-cofinite
modules, a $\mathcal{Z}$-cofinite module can have infinitely many associated prime ideals.

\begin{example} \label{3.3}
\begin{enumerate}
\item[(i)] Let $R$ be a ring with infinitely many maximal ideals. Let $\mathcal{Z}:= \Max(R)$, and $M:= \bigoplus_{\mathfrak{m}
\in \mathcal{Z}} R/\mathfrak{m}$. Then $F(\mathcal{Z})= \{\mathfrak{a} \lhd R \suchthat \dim(R/\mathfrak{a})\leq 0 \}$. Thus for
any given proper ideal $\mathfrak{a} \in F(\mathcal{Z})$, there are finitely many maximal ideals in $V(\mathfrak{a})$, say
$\mathfrak{m}_{1},..., \mathfrak{m}_{n}$. It follows that $\Gamma_{\mathfrak{a}}(M)=\bigoplus_{i=1}^{n} R/\mathfrak{m}_{i}$.
Now $\Gamma_{\mathfrak{a}}(M)$ is artinian, while $\Gamma_{\mathcal{Z}}(M)=M$ fails to be artinian as it contains infinitely
many direct summands.

\item[(ii)] Let $R$ be a Gorenstein ring of finite dimension $d$ such that $\mathcal{Z}:= \{\mathfrak{m}\in \Max(R)
\suchthat \Ht(\mathfrak{m})=d \}$ is an infinite set. Then by \cite[Remark 2.12]{HD}, it turns out that
\[
    H_{\mathcal{Z}}^{j}(R)\cong
\begin{dcases}
    \bigoplus_{\mathfrak{m}\in \mathcal{Z}}E_{R}(R/\mathfrak{m}) & \text{if } j=d\\
    0              & \text{if } j\neq d.
\end{dcases}
\]
Thus $\Ext^{i}_{R}\left(R/\mathfrak{a},H^{j}_{\mathcal{Z}}(R)\right)$ is finitely generated for every $\mathfrak{a} \in F(\mathcal{Z})$ and every $i,j \geq 0$. It follows that $H^{j}_{\mathcal{Z}}(R)$ is $\mathcal{Z}$-cofinite for every $j \geq 0$, whereas $H^{d}_{\mathcal{Z}}(R)$ is not artinian and $\Ass_R\left(H^{d}_{\mathcal{Z}}(R)\right)$ is not finite.
\end{enumerate}
\end{example}

\begin{lemma} \label{3.4}
Let $R$ be a semilocal ring with Jacobson radical $\mathfrak{J}$, $\mathcal{Z}$ a stable under specialization subset of $\Spec(R)$, and $M$ an $R$-module. Let $(-)^{\vee}:= \Hom_{R}\left(-,E_{R}\left(R/\mathfrak{J}\right)\right)$ be the Matlis duality functor. Then the following assertions are equivalent for any given $r \geq 0$:
\begin{enumerate}
\item[(i)] $H^{i}_{\mathcal{Z}}\left(M^{\vee}\right)$ is an artinian $R$-module for every $0 \leq i \leq r$.
\item[(ii)] $\Ext^{i}_{R}\left(R/\mathfrak{a},M^{\vee}\right)$ is an artinian $R$-module for every $\mathfrak{a} \in F(\mathcal{Z})$ and every $0 \leq i \leq r$.
\item[(iii)] $\Tor^{R}_{i}(R/\mathfrak{a},M)$ is a finitely generated $R$-module for every $\mathfrak{a} \in F(\mathcal{Z})$ and every $0 \leq i \leq r$.
\end{enumerate}
\end{lemma}

\begin{prf}
(i) $\Leftrightarrow$ (ii): Follows from Lemma \ref{3.2}.

One can easily see that $E_{R}\left(R/\mathfrak{J}\right) \cong \bigoplus_{\mathfrak{m}\in \Max(R)}E_{R}(R/\mathfrak{m})$
is an artinian injective cogenerator for $R$. Fix $\mathfrak{a} \in F(\mathcal{Z})$ and $0 \leq i \leq r$, and let $N:=
\Tor^{R}_{i}(R/\mathfrak{a},M)$ for the rest of the proof. Then by \cite[Corollary 10.63]{Ro}, we have $N^{\vee}\cong \Ext^{i}_{R}
\left(R/\mathfrak{a},M^{\vee}\right)$.

(ii) $\Rightarrow$ (iii): Let $\Max(R)=\{\mathfrak{m}_1,\dots ,\mathfrak{m}_n\}$ and set $T:=\widehat{R}^{\mathfrak{J}}$
and $T_j:=\widehat{R_{\mathfrak{m}_{j}}}$ for every $j=1,\dots, n$. We know that $T\cong \prod_{j=1}^nT_j$ and $T$ is a
$\mathfrak{J}T$-adically complete semilocal ring with $\Max(T)=\{\mathfrak{m}_{j}T \suchthat j=1,\dots, n\}$.
Any $\mathfrak{J}$-torsion $R$-module possesses a $T$-module structure in such a way that a subset is an $R$-submodule if and only
if it is a $T$-submodule. In particular, $N^{\vee}$ and $E_{R}\left(R/\mathfrak{J}\right)$ are artinian $T$-modules.
Moreover, one may easily check that the two $T$-modules $E_{R}\left(R/\mathfrak{J}\right)$ and $E_{T}\left(T/\mathfrak{J}T\right)$
are isomorphic and $\mathfrak{J}T$ is the Jacobson radical of $T$. Putting everything together, we obtain:
\begin{equation*}
\begin{split}
N^{\vee} & \cong \Hom_{R}\left(N,E_{R}\left(R/\mathfrak{J}\right)\right) \\
 & \cong \Hom_{R}\left(N,E_{T}\left(T/\mathfrak{J}T\right)\right) \\
 & \cong \Hom_{R}\left(N,\Hom_{T}\left(T,E_{T}\left(T/\mathfrak{J}T\right)\right)\right)  \\
 & \cong \Hom_{T}\left(N\otimes_{R} T,E_{T}\left(T/\mathfrak{J}T\right)\right). \\
 \end{split}
\end{equation*}
Applying the Matlis Duality Theorem over the ring $T$ \cite[Proposition 4 (c)]{CW}, we deduce that $N\otimes_{R}T$ is a
finitely generated $T$-module and by the faithfully flatness of the completion map $\theta^{\mathfrak{J}}_{R}: R \rightarrow T$,
we infer that $N$ is a finitely generated  $R$-module.

(iii) $\Rightarrow$ (ii): There is an exact sequence $$R^{n}\rightarrow N \rightarrow 0,$$ which yields the exact sequence
$$0 \rightarrow N^{\vee} \rightarrow E_{R}\left(R/\mathfrak{J}\right)^{n}.$$ It follows that $\Ext^{i}_{R}\left(R/\mathfrak{a},
M^{\vee}\right) \cong N^{\vee}$ is an artinian $R$-module.
\end{prf}

For any given $R$-module $N$, one may easily check that $\cd(\mathcal{Z},N) \leq \dim(R)$.

\begin{corollary} \label{3.5}
Let $R$ be a semilocal ring, $\mathcal{Z}$ a stable under specialization subset of $\Spec(R)$, and $M$ an $R$-module with
$\Supp_{R}(M) \subseteq \mathcal{Z}$. Then the following assertions are equivalent:
\begin{enumerate}
\item[(i)] $M$ is $\mathcal{Z}$-cofinite.
\item[(ii)] $\Tor^{R}_{i}(R/\mathfrak{a},M)$ is finitely generated for every $\mathfrak{a} \in F(\mathcal{Z})$ and every $i \geq 0$.
\item[(iii)] $H^{i}_{\mathcal{Z}}\left(M^{\vee}\right)$ is artinian for every $i \geq 0$.
\item[(iv)] $H^{i}_{\mathcal{Z}}\left(M^{\vee}\right)$ is artinian for every $0 \leq i \leq \cd(\mathcal{Z},M^{\vee})$.
\item[(v)] $\Ext^{i}_{R}\left(R/\mathfrak{a},M^{\vee}\right)$ is artinian for every $\mathfrak{a} \in F(\mathcal{Z})$ and every $0 \leq i \leq \cd(\mathcal{Z},M^{\vee})$.
\item[(vi)] $\Tor^{R}_{i}(R/\mathfrak{a},M)$ is finitely generated for every $\mathfrak{a} \in F(\mathcal{Z})$ and every $0 \leq i \leq \cd(\mathcal{Z},M^{\vee})$.
\end{enumerate}
\end{corollary}

\begin{prf}
(i) $\Leftrightarrow$ (ii): Follows from Lemma \ref{2.2} (iii).

(ii) $\Leftrightarrow$ (iii): Follows from Lemma \ref{3.4}.

(iii) $\Leftrightarrow$ (iv): Obvious.

(iv) $\Leftrightarrow$ (v) $\Leftrightarrow$ (vi): Follows from Lemma \ref{3.4}.
\end{prf}

The following two lemmas are our essential tools in the proof of the second part of Theorem \ref{3.11}.

\begin{lemma} \label{3.6}
Let $\mathcal{Z}$ be a stable under specialization subset of $\Spec(R)$, and $M$ an $R$-module with $\Supp_{R}(M) \subseteq \mathcal{Z}$ and $\dim_R(M)\leq 1$.
Then the following assertions are equivalent:
\begin{enumerate}
\item[(i)] $M$ is $\mathcal{Z}$-cofnite.
\item[(ii)] $\Hom_{R}(R/\mathfrak{a},M)$ and $\Ext^{1}_{R}(R/\mathfrak{a},M)$ are finitely generated for every $\mathfrak{a} \in F(\mathcal{Z})$.
\end{enumerate}
\end{lemma}

\begin{prf} Follows from the equivalence (i)$\Leftrightarrow$(iii) in \cite[Lemma 2.6]{AB}.
\end{prf}

In the sequel, we use the straightforward observation that if any two modules in a short exact sequence are $\mathcal{Z}$-cofinite, then so is the third.

\begin{lemma} \label{4.7}
Let $\mathcal{Z}$ be a stable under specialization subset of $\Spec(R)$. Assume that either
\begin{enumerate}
\item[(i)] An $R$-module $M$ is $\mathcal{Z}$-cofinite whenever $\Supp_R(M)\subseteq \mathcal{Z}$ and $\Ext_R^i(R/\fa,M)$ is finitely generated for $i=0,1$ and every
$\mathfrak{a}\in F(\mathcal{Z})$; or
\item[(ii)] An $R$-module $M$ is $\mathcal{Z}$-cofinite whenever $\Supp_R(M)\subseteq \mathcal{Z}$ and $\Tor^R_i(R/\fa,M)$ is finitely generated for $i=0,1$ and every
$\mathfrak{a}\in F(\mathcal{Z}),$
\end{enumerate}
holds. Then $\mathcal{M}(R,\mathcal{Z})_{cof}$ is an abelian subcategory of $\mathcal{M}(R)$.
\end{lemma}

\begin{prf}
(i): Let $f:M \rightarrow N$ be an $R$-homomorphism between $\mathcal{Z}$-cofinite modules. We should prove that both $\ker f$ and $\coker f$ are $\mathcal{Z}$-cofinite.
In view of the short exact sequences
\begin{equation} \label{eq:4.7.1}
0 \rightarrow \im f \rightarrow N \rightarrow \coker f \rightarrow 0,
\end{equation}
and
\begin{equation} \label{eq:4.7.2}
0 \rightarrow \ker f \rightarrow M \rightarrow \im f \rightarrow 0,
\end{equation}
it suffices to show that $\ker f$ is $\mathcal{Z}$-cofinite. Let $\mathfrak{b}\in F(\mathcal{Z}).$ From \eqref{eq:4.7.1}, we deduce that $\Hom_{R}(R/\mathfrak{b},\im f)$
is finitely generated. Now, \eqref{eq:4.7.2} yields the exact sequence $$0 \rightarrow \Hom_{R}(R/\mathfrak{b},\ker f) \rightarrow \Hom_{R}(R/\mathfrak{b},M) \rightarrow \Hom_{R}(R/\mathfrak{b},\im f) \rightarrow $$$$ \Ext_{R}^{1}(R/\mathfrak{b},\ker f) \rightarrow \Ext_{R}^{1}(R/\mathfrak{b},M).$$ Thus $\Hom_{R}(R/\mathfrak{b},\ker f)$ and $\Ext_{R}^{1}(R/\mathfrak{b},\ker f)$ are finitely generated. Next, our assumption in (i) implies that $\ker f$ is $\mathcal{Z}$-cofinite. Note that  $\Supp_R(\ker f)\subseteq
\mathcal{Z}$.

(ii): It is similar to the proof of (i), and so we leave it to the reader.
\end{prf}

To prove the third part of Theorem \ref{3.11}, we need Lemmas \ref{3.7}, \ref{3.8}, \ref{3.9}, and \ref{3.10}.

\begin{lemma} \label{3.7}
Let $S$ be a module-finite $R$-algebra. Let $\mathfrak{a}$ be an ideal of $R$, and $M$ an $S$-module. Then the $R$-module $\Ext^{i}_{R}(R/\mathfrak{a},M)$ is finitely
generated for every $i\geq 0$ if and only if the $S$-module $\Ext^{i}_{S}(S/\mathfrak{a}S,M)$ is finitely generated for every $i\geq 0$.
\end{lemma}

\begin{prf}
The proof of \cite[Proposition 2]{DM} establishes the claim. Note that the assumption on the supports is not used in that proof.
\end{prf}

Let $\fa$ be an ideal of $R$. An $R$-module $M$ is said to be $\mathfrak{a}$-Ext-finite if $\Ext^{i}_{R}(R/\mathfrak{a},M)$ is finitely generated for every $i\geq 0$.

\begin{lemma} \label{3.8}
Let $\frak{a}$ be a proper ideal of $R$, and $M$ an $R$-module. Suppose that $\dim(R) \leq 1$. Then the following assertions hold:
\begin{enumerate}
\item[(i)] If $(0:_{M}\mathfrak{a})$ is finitely generated, then $M$ is $\mathfrak{a}$-Ext-finite.
\item[(ii)] The class of $\mathfrak{a}$-Ext-finite $R$-modules is closed under taking submodules, quotients, and extensions.
\end{enumerate}
\end{lemma}

\begin{prf}
(i): There is an integer $n \geq 1$ such that $\Gamma_{\mathfrak{a}}(R)=(0:_{R}\mathfrak{a}^{n})$. Then $\bar{M}:= M/(0:_{M}\mathfrak{a}^{n})$ is a module over the ring
$\bar{R}:= R/\Gamma_{\mathfrak{a}}(R)$. Let $\bar{\mathfrak{a}}$ be the image of $\mathfrak{a}$ in $\bar{R}$. Then $\bar{\mathfrak{a}}$ contains an $\bar{R}$-regular
element and thus $\dim(\bar{R}/\bar{\mathfrak{a}})=0$. We note that as $(0:_{M}\mathfrak{a})$ is finitely generated, one may check that $(0:_{M}\mathfrak{a}^{i})$ is
finitely generated for every $i\geq 1$. So, $(0:_{\bar{M}}\bar{\mathfrak{a}})=(0:_{M}\mathfrak{a}^{n+1})/(0:_{M}\mathfrak{a}^{n})$ is a finitely generated $R$-module.
It follows that $(0:_{\bar{M}}\bar{\mathfrak{a}})$
is a finitely generated $\bar{R}/\bar{\mathfrak{a}}$-module. As $\bar{R}/\bar{\mathfrak{a}}$ is artinian, we see that $(0:_{\bar{M}}\bar{\mathfrak{a}})$ is artinian as an $\bar{R}/\bar{\mathfrak{a}}$-module, and thus as an $\bar{R}$-module. By \cite[Theorem 1.3]{Me3}, it follows that $\Gamma_{\bar{\mathfrak{a}}}(\bar{M})$ is an artinian
$\bar{R}$-module. But $\dim_{\bar{R}}(\bar{M}) \leq 1$, so \cite[Proposition 5.1]{Me2} implies that $H^{1}_{\bar{\mathfrak{a}}}(\bar{M})$ is an artinian
$\bar{\mathfrak{a}}$-cofinite $\bar{R}$-module. It now follows from Lemma \ref{3.2} that $\Ext^{i}_{\bar{R}}(\bar{R}/\bar{\mathfrak{a}},\bar{M})$ is an artinian
$\bar{R}$-module for every $i \geq 0$. But $\bar{R}/\bar{\mathfrak{a}}$ is artinian, so $\Ext^{i}_{\bar{R}}(\bar{R}/\bar{\mathfrak{a}},\bar{M})$ has finite length for every
$i \geq 0$. It follows from Lemma \ref{3.7} that $\Ext^{i}_{R}(R/\mathfrak{a},\bar{M})$ is a finitely generated $R$-module for every $i \geq 0$. The short exact sequence
$$0 \rightarrow (0:_{M}\mathfrak{a}^{n}) \rightarrow M \rightarrow \bar{M} \rightarrow 0$$
implies that $\Ext^{i}_{R}(R/\mathfrak{a},M)$ is a finitely generated $R$-module for every $i \geq 0$.

(ii): Follows from (i).
\end{prf}

\begin{lemma} \label{3.9}
Let $\mathfrak{a}$ be an ideal of $R$ and $a \in \mathfrak{a}$. If $L$ is an $R$-module such that $L/aL$ and $(0:_{L}a)$ are $\mathfrak{a}$-Ext-finite, then $L$ is
$\mathfrak{a}$-Ext-finite.
\end{lemma}

\begin{prf}
Apply \cite[Corollary 3.3]{Me2} to the $R$-homomorphism $f=a\text{Id}_L$. We deduce that $L$ is $\mathfrak{a}$-Ext-finite. Note that $\Ext^{i}_{R}(R/\mathfrak{a},f)=0$
for every $i\geq 0$ as $a \in \mathfrak{a}$.
\end{prf}

\begin{lemma} \label{3.10}
Let $\frak{a}$ be an ideal of $R$. Suppose that $\dim(R) \leq 2$, and there exists an $a \in \mathfrak{a}$ with $\dim(R/aR) \leq 1$. If $M$ is an $\mathfrak{a}$-Ext-finite
$R$-module, then $M/aM$ and $(0:_{M}a)$ are $\mathfrak{a}$-Ext-finite.
\end{lemma}

\begin{prf}
Take elements $a_{1},...,a_{n} \in \mathfrak{a}$ such that $a_{1}=a$ and $\mathfrak{a}=(a_{1},a_{2},...,a_{n})$. Then by \cite[Theorem 2.1]{Me2}, $H^{i}(a_{1},...,a_{n};M)$ is
finitely generated for every $i \geq 0$, so it is $\mathfrak{a}$-Ext-finite for every $i \geq 0$.

Let $i \geq 0$ and $L := H^{i}(a_{1},...,a_{n-1};M)$. Let $\bar{\mathfrak{a}}$ be the image of $\mathfrak{a}$ in $\bar{R}:=R/aR$. By Lemma \ref{3.7}, $H^{i}(a_{1},...,a_{n};M)$ is $\bar{\mathfrak{a}}$-Ext-finite for every $i \geq 0$. In the exact sequence
$$H^{i}(a_{1},...,a_{n};M) \rightarrow H^{i}(a_{1},...,a_{n-1};M) \xrightarrow{a_{n}} H^{i}(a_{1},...,a_{n-1};M) \rightarrow H^{i+1}(a_{1},...,a_{n};M),$$
the outer terms are $\bar{\mathfrak{a}}$-Ext-finite. By Lemma \ref{3.8} (ii), $L/a_{n}L$ and $(0:_{L}a_{n})$ are $\bar{\mathfrak{a}}$-Ext-finite. Hence,
Lemma \ref{3.7} implies that $L/a_{n}L$ and $(0:_{L}a_{n})$ are $\mathfrak{a}$-Ext-finite. By Lemma \ref{3.9}, $L$ is $\mathfrak{a}$-Ext-finite.

Continuing in this fashion, we infer that $H^{i}(a_{1};M)$ is $\mathfrak{a}$-Ext-finite for every $i \geq 0$, and so $M/aM$ and $(0:_{M}a)$ are $\mathfrak{a}$-Ext-finite.
\end{prf}

\begin{remark} \label{3.10.1}
Let $\fa$ be a proper ideal of $R$ and $f:M \rightarrow N$ an $R$-homomorphism between $\mathfrak{a}$-Ext-finite modules. We show that in case we would like to
prove that $\ker f$ and $\coker f$ are $\mathfrak{a}$-Ext-finite, we may additionally assume that $\mathfrak{a}$ contains an $R$-regular element. Indeed, there
is an integer $n \geq 1$ such that $\Gamma_{\mathfrak{a}}(R) = (0:_{R}\mathfrak{a}^{n})$. Let $\bar{R} = R/ (0:_{R}\mathfrak{a}^{n})$ and $\bar{\mathfrak{a}}=
\mathfrak{a}\bar{R}$. It is clear that $\depth_{\bar{R}}(\bar{\mathfrak{a}},\bar{R})>0$. Consider the following commutative diagram with exact rows:
\[\begin{tikzpicture}[every node/.style={midway},]
  \matrix[column sep={3em}, row sep={3em}]
  {\node(1) {$0$}; & \node(2) {$(0:_{M}\mathfrak{a}^{n})$}; & \node(3) {$M$}; & \node(4) {$M/(0:_{M}\mathfrak{a}^{n})$}; & \node(5) {$0$};\\
  \node(6) {$0$}; & \node(7) {$(0:_{N}\mathfrak{a}^{n})$}; & \node(8) {$N$}; & \node(9) {$N/(0:_{N}\mathfrak{a}^{n})$}; & \node(10) {$0,$};\\};
  \draw[decoration={markings,mark=at position 1 with {\arrow[scale=1.5]{>}}},postaction={decorate},shorten >=0.5pt] (2) -- (7) node[anchor=west] {$\tilde{f}$};
  \draw[decoration={markings,mark=at position 1 with {\arrow[scale=1.5]{>}}},postaction={decorate},shorten >=0.5pt] (3) -- (8) node[anchor=west] {$f$};
  \draw[decoration={markings,mark=at position 1 with {\arrow[scale=1.5]{>}}},postaction={decorate},shorten >=0.5pt] (4) -- (9) node[anchor=west] {$\bar{f}$};
  \draw[decoration={markings,mark=at position 1 with {\arrow[scale=1.5]{>}}},postaction={decorate},shorten >=0.5pt] (1) -- (2) node[anchor=south] {};
  \draw[decoration={markings,mark=at position 1 with {\arrow[scale=1.5]{>}}},postaction={decorate},shorten >=0.5pt] (2) -- (3) node[anchor=south] {};
  \draw[decoration={markings,mark=at position 1 with {\arrow[scale=1.5]{>}}},postaction={decorate},shorten >=0.5pt] (3) -- (4) node[anchor=south] {};
  \draw[decoration={markings,mark=at position 1 with {\arrow[scale=1.5]{>}}},postaction={decorate},shorten >=0.5pt] (4) -- (5) node[anchor=south] {};
  \draw[decoration={markings,mark=at position 1 with {\arrow[scale=1.5]{>}}},postaction={decorate},shorten >=0.5pt] (6) -- (7) node[anchor=south] {};
  \draw[decoration={markings,mark=at position 1 with {\arrow[scale=1.5]{>}}},postaction={decorate},shorten >=0.5pt] (7) -- (8) node[anchor=south] {};
  \draw[decoration={markings,mark=at position 1 with {\arrow[scale=1.5]{>}}},postaction={decorate},shorten >=0.5pt] (8) -- (9) node[anchor=south] {};
  \draw[decoration={markings,mark=at position 1 with {\arrow[scale=1.5]{>}}},postaction={decorate},shorten >=0.5pt](9) -- (10) node[anchor=south] {};
\end{tikzpicture}\]
where $\tilde{f}$ and $\bar{f}$ are induced by $f$ in the obvious way. Since the $R$-modules $(0:_{M}\mathfrak{a})$ and $(0:_{N}\mathfrak{a})$ are
finitely generated, it can be seen that the $R$-modules $(0:_{M}\mathfrak{a}^{n})$ and $(0:_{M}\mathfrak{a}^{n})$ are also finitely generated.
So, from the rows of the above commutative diagram, one deduces that the $R$-modues $M/(0:_{M}\mathfrak{a}^{n})$ and $N/(0:_{N}\mathfrak{a}^{n})$ are
$\mathfrak{a}$-Ext-finite.

Applying the Snake Lemma to the above diagram, we get the exact sequence
$$0 \rightarrow \ker \tilde{f} \rightarrow \ker f \rightarrow \ker \bar{f} \rightarrow \coker \tilde{f} \rightarrow \coker f \rightarrow \coker \bar{f} \rightarrow 0.$$
Now, $\ker \tilde{f}$ and $\coker \tilde{f}$ are finitely generated $R$-modules. Hence, $\ker f$ and $\coker f$ are $\mathfrak{a}$-Ext-finite if and only if
$\ker \bar{f}$ and $\coker \bar{f}$ are $\mathfrak{a}$-Ext-finite. But, Lemma \ref{3.7} yields that $\ker \bar{f}$ and $\coker \bar{f}$ are $\mathfrak{a}$-Ext-finite if
and only if they are $\mathfrak{\bar{a}}$-Ext-finite.
\end{remark}

Finally, we are ready to prove the main result of this section.

\begin{theorem} \label{3.11}
Let $\mathcal{Z}$ be a stable under specialization subset of $\Spec (R)$. Assume that either
\begin{enumerate}
\item[(i)]  $R$ is semilocal with $\cd(\mathcal{Z},R) \leq 1$, or
\item[(ii)] $\dim(\mathcal{Z}) \leq 1$, or
\item[(iii)] $\dim(R) \leq 2$.
\end{enumerate}
Then $\mathcal{M}(R,\mathcal{Z})_{cof}$ is an abelian subcategory of $\mathcal{M}(R)$.
\end{theorem}

\begin{prf}
(i): Let $M$ be an $R$-module, and $N$ a finitely generated $R$-module such that $\Supp_{R}(M)\subseteq \Supp_{R}(N)$. We claim that $\cd(\mathcal{Z},M)\leq \cd(\mathcal{Z},N)$. Since $H^{i}_{\mathcal{Z}}(-)$ commutes with direct limits and $M$ can be written as a direct limit of its finitely generated submodules, we may assume that $M$ is finitely generated. Now, the proof is a straightforward adaptation of the argument given in \cite[Theorem 2.2]{DNT}. In particular, $\cd(\mathcal{Z},L)\leq \cd(\mathcal{Z},R)\leq 1$ for every $R$-module $L$.
Thus, the assertion is immediate by Corollary \ref{3.5} and Lemma \ref{4.7} (ii).

(ii): Since  $\dim(\mathcal{Z})\leq 1$, it turns out that $\dim_R(M)\leq 1$ for every $\mathcal{Z}$-cofinite $R$-module $M$. Hence, the claim follows from Lemma \ref{3.6} and
Lemma \ref{4.7} (i).

(iii): Let $f:M \rightarrow N$ be an $R$-homomorphism between $\mathcal{Z}$-cofinite $R$-modules. We need to show that both $\ker f$ and $\coker f$ are $\mathfrak{a}$-Ext-finite for
every $\mathfrak{a} \in F(\mathcal{Z})$. Fix $\mathfrak{a} \in F(\mathcal{Z})$. In light of Remark \ref{3.10.1}, we may assume that there is an $R$-regular element $a \in \mathfrak{a}$. Besides, in view of the short exact sequences
$$0\lo \ker f\lo M\lo \im f \lo 0$$
and
$$0\lo \im f\lo N\lo \coker f \lo 0,$$
it is sufficient to show that $\im f$ is $\mathfrak{a}$-Ext-finite. The $R$-modules $M/aM$ and $(0:_{N}a)$ are $\mathfrak{a}$-Ext-finite by Lemma \ref{3.10}. Hence Lemma \ref{3.7} implies that $M/aM$ and $(0:_{N}a)$ are $\mathfrak{a}/aR$-Ext-finite over the ring $\bar{R}:=R/aR$ which has dimension at most one. We note that $\im f / a\im f$ is a homomorphic image of $M/aM$, and $(0:_{\im f}a)$ is a submodule of $(0:_{N}a)$. Therefore, $\im f / a\im f$ and $(0:_{\im f}a)$ are $\mathfrak{a}/aR$-Ext-finite by Lemma \ref{3.8} (ii). Hence,
they are $\mathfrak{a}$-Ext-finite by Lemma \ref{3.7}. Now, Lemma \ref{3.9} implies that $\im f$ is $\mathfrak{a}$-Ext-finite.
\end{prf}

In Theorem \ref{3.11} (i), the assumption that $R$ is semilocal is somehow not desirable. Accordingly, we pose the following question.

\begin{question} \label{3.12}
Let $\mathcal{Z}$ be a stable under specialization subset of $\Spec(R)$ such that $\cd(\mathcal{Z},R) \leq 1$. Is $\mathcal{M}(R,\mathcal{Z})_{cof}$ an abelian subcategory of $\mathcal{M}(R)$?
\end{question}

\section{Proof of Theorem 1.6}

In this section, we intend to prove Theorem \ref{1.6}; see Theorem \ref{4.8}. To this end, Lemma \ref{4.6} is our main tool. In order to apply it in the situation of
Theorem \ref{4.8}, one has to use Theorem \ref{3.11}, and Lemmas \ref{4.3} and \ref{4.7}. For proving Lemma \ref{4.3}, we need Lemmas \ref{4.1} and \ref{4.2}.

Given a stable under specialization subset $\mathcal{Z}$ of $\Spec(R)$ and a finitely generated $R$-module $M$, we remind that
$$\depth_{R}(\mathcal{Z},M):= \inf \left\{\depth_{R}(\mathfrak{a},M) \suchthat \mathfrak{a} \in F(\mathcal{Z}) \right\}.$$

\begin{lemma} \label{4.1}
Let $\mathcal{Z}$ be a stable under specialization subset of $\Spec(R)$, and $M$ a finitely generated $R$-module. Then
$\Supp_{R}(M)\cap \mathcal{Z}=\bigcup_{i=0}^{\infty} \Supp_{R}\left(H^{i}_{\mathcal{Z}}(M)\right)$.
\end{lemma}

\begin{prf}
Let $\mathfrak{p} \in \Supp_{R}(M)\cap \mathcal{Z}$. It is straightforward to see that the set
$$\mathcal{Z}_{\mathfrak{p}}:= \left\{\mathfrak{q}R_{\mathfrak{p}} \suchthat \mathfrak{q}\in \mathcal{Z} \text{ and } \mathfrak{q}\subseteq \mathfrak{p}\right\}$$
is a stable under specialization subset of $\Spec(R_{\frak p})$. It is clear that $\mathfrak{p}R_{\mathfrak{p}} \in \mathcal{Z}_{\mathfrak{p}} \cap \Supp_{R_{\mathfrak{p}}}\left(M_{\mathfrak{p}}\right)$, so $\depth_{R_{\mathfrak{p}}}\left(\mathfrak{p}R_{\mathfrak{p}},M_{\mathfrak{p}}\right) < \infty$, and thus $s:= \depth_{R_{\mathfrak{p}}}\left(\mathcal{Z}_{\mathfrak{p}},M_{\mathfrak{p}}\right) < \infty$. However by \cite[Proposition 5.5]{Bi},
$$\depth_{R_{\mathfrak{p}}}\left(\mathcal{Z}_{\mathfrak{p}},M_{\mathfrak{p}}\right)= \inf \left\{i \in \mathbb{Z} \suchthat H^{i}_{\mathcal{Z_{\mathfrak{p}}}}\left(M_{\mathfrak{p}}\right) \neq 0 \right\},$$ so $H^{s}_{\mathcal{Z_{\mathfrak{p}}}}\left(M_{\mathfrak{p}}\right) \neq 0$. One may check that $H^{s}_{\mathcal{Z}}(M)_{\mathfrak{p}} \cong H^{s}_{\mathcal{Z_{\mathfrak{p}}}}\left(M_{\mathfrak{p}}\right)$, and so $\mathfrak{p}\in \bigcup_{i=0}^{\infty} \Supp_{R}\left(H^{i}_{\mathcal{Z}}(M)\right)$. The reverse inclusion is immediate.
\end{prf}

\begin{lemma} \label{4.2}
Let $\mathcal{Z}$ be a stable under specialization subset of $\Spec(R)$, and $M$ a finitely generated $R$-module. Then $\dim\left(\Supp_{R}\left(H^{i}_{\mathcal{Z}}(M)\right)\right) \leq \dim_{R}(M)-1$ for every $i \geq 1$.
\end{lemma}

\begin{prf}
Clearly, we may assume that $\dim_{R}(M) < \infty$. Since $$H_{\mathcal{Z}}^i(M) \cong \underset{\fa\in F(\mathcal{Z})}{\varinjlim}H_{\mathfrak{a}}^i(M)$$ for every $i \geq 0$, it suffices to show that $\dim \left(\Supp_{R}\left(H_{\mathfrak{a}}^i(M)\right)\right) \leq \dim_{R}(M)-1$ for every $\fa\in F(\mathcal{Z})$ and every $i \geq 1$. As $H^{i}_{\mathfrak{a}}(M) \cong H^{i}_{\mathfrak{a}}\left(M/\Gamma_{\mathfrak{a}}(M)\right)$ for every $i \geq 1$, we may assume that $\Gamma_{\mathfrak{a}}(M)=0$. Consequently, we conclude that $\mathfrak{a}$ contains a nonzerodivisor $r$ on $M$. As $\Supp_{R}\left(H^{i}_{\mathfrak{a}}(M)\right) \subseteq \Supp_{R}(M/\mathfrak{a}M)$, we have
\begin{equation*}
\begin{split}
\dim\left(\Supp_{R}\left(H^{i}_{\mathfrak{a}}(M)\right)\right) & \leq \dim_{R}(M/\mathfrak{a}M) \\
 & \leq \dim_{R}(M/rM) \\
 & \leq \dim_{R}(M)-1
\end{split}
\end{equation*}
for every $i \geq 1$.
\end{prf}

\begin{lemma} \label{4.3}
Let $\mathcal{Z}$ be a stable under specialization subset of $\Spec(R)$, and $M$ a finitely generated $R$-module.
Assume that either
\begin{enumerate}
\item[(i)]  $\cd(\mathcal{Z},R) \leq 1$, or
\item[(ii)] $\dim(\mathcal{Z}) \leq 1$, or
\item[(iii)] $\dim\left(\Supp_{R}\left(H^{i}_{\mathcal{Z}}(M)\right)\right) \leq 1$ for every $i\geq 0$, or
\item[(iv)] $\dim_R(M)\leq 2$.
\end{enumerate}
Then $H^{i}_{\mathcal{Z}}(M)$ is $\mathcal{Z}$-cofinite for every $i \geq 0$.
\end{lemma}

\begin{prf}
(i): See \cite[Proposition 2.7]{MS}.

(ii): Clearly, $\Supp_{R}\left(H^{i}_{\mathcal{Z}}(M)\right) \subseteq \mathcal{Z}$ for every $i \geq 0$. Thus, it remains to show that $\Ext^{j}_{R}\left(R/\mathfrak{a},H^{i}_{\mathcal{Z}}(M)\right)$ is finitely generated for every $\mathfrak{a} \in F(\mathcal{Z})$ and every $i,j \geq 0$.

Fix $\mathfrak{b}\in F(\mathcal{Z})$. By induction on $i$, for a given $R$-module $N$ with $\Ext^{j}_{R}\left(R/\mathfrak{b},N\right)$ finitely generated
for every $j\geq 0$, we show that $\Ext^{j}_{R}\left(R/\mathfrak{b},H^{i}_{\mathcal{Z}}(N)\right)$ is finitely generated for every $j\geq 0$. The short
exact sequence $$0 \rightarrow \Gamma_{\mathcal{Z}}(N) \rightarrow N \rightarrow M/\Gamma_{\mathcal{Z}}(N) \rightarrow 0,$$
yields the exact sequence
$$0 \rightarrow \Hom_{R}\left(R/\mathfrak{b},\Gamma_{\mathcal{Z}}(N)\right) \rightarrow \Hom_{R}(R/\mathfrak{b},N) \rightarrow \Hom_{R}\left(R/\mathfrak{b},N/\Gamma_{\mathcal{Z}}(N)\right) \rightarrow $$$$ \Ext_{R}^{1}\left(R/\mathfrak{b},\Gamma_{\mathcal{Z}}(N)\right) \rightarrow \Ext_{R}^{1}(R/\mathfrak{b},N).$$
It can be seen by inspection that $\Hom_{R}\left(R/\mathfrak{b},N/\Gamma_{\mathcal{Z}}(N)\right)=0$. Hence the above exact sequence shows that the $R$-modules $\Hom_{R}\left(R/\mathfrak{b},\Gamma_{\mathcal{Z}}(N)\right)$ and $\Ext_{R}^{1}\left(R/\mathfrak{b},\Gamma_{\mathcal{Z}}(N)\right)$ are finitely generated. Therefore, by Lemma \ref{3.6} the case $i=0$ holds true.

Now, suppose that $i\geq 1$ and make the obvious induction hypothesis. From the exact sequence
$$\Ext_{R}^{j}(R/\mathfrak{b},N) \rightarrow \Ext_{R}^{j}\left(R/\mathfrak{b},N/\Gamma_{\mathcal{Z}}(N)\right) \rightarrow \Ext_{R}^{j+1}\left(R/\mathfrak{b},\Gamma_{\mathcal{Z}}(N)\right),$$
using the base case $i=0$, we deduce that $\Ext_{R}^{j}\left(R/\mathfrak{b},N/\Gamma_{\mathcal{Z}}(N)\right)$ is finitely generated for every $j\geq 0$. Since $H^{i}_{\mathcal{Z}}(N)\cong H^{i}_{\mathcal{Z}}\left(N/\Gamma_{\mathcal{Z}}(N)\right)$ for every $i \geq 1$, we may assume that $\Gamma_{\mathcal{Z}}(N)=0$. Let $E:= E_{R}(N)$ and $L:=E/N$. We have $\Gamma_{\mathcal{Z}}(E)\cong E_{R}\left(\Gamma_{\mathcal{Z}}(N)\right)=0$, and $\Hom_{R}(R/\mathfrak{b},E)=0$. Then from the short exact sequence
$$0 \rightarrow N \rightarrow E \rightarrow L \rightarrow 0,$$
we conclude that $H^{k}_{\mathcal{Z}}(N) \cong H^{k-1}_{\mathcal{Z}}(L)$ and $\Ext^{k}_{R}(R/\mathfrak{b},N) \cong \Ext^{k-1}_{R}(R/\mathfrak{b},L)$ for every $k\geq 1$.
Hence the assumption is satisfied by $L$, and thus $\Ext^{j}_{R}\left(R/\mathfrak{b},H^{i-1}_{\mathcal{Z}}(L)\right)$ is finitely generated for every $j\geq 0$.

(iii): By the assumption and Lemma \ref{4.1}, we have $\dim\left(\Supp_{R}(M)\cap \mathcal{Z}\right) \leq 1$. Set
$$\widetilde{\mathcal{Z}}:= \left\{\frac{\mathfrak{p}}{\ann_{R}(M)} \suchthat \mathfrak{p} \in \Supp_{R}(M)\cap
\mathcal{Z} \right\},$$ and $S:= R/ \ann_{R}(M)$. Then, it is straightforward to see that $\widetilde{\mathcal{Z}}$
is a stable under specialization subset of $\Spec(S)$ with $\dim\left(\widetilde{\mathcal{Z}}\right) \leq 1$, and
$F\left(\widetilde{\mathcal{Z}}\right)=\left\{\mathfrak{a}S \suchthat \mathfrak{a} \in F(\mathcal{Z}) \right\}$.
In addition, we have
\begin{equation*}
\begin{split}
H^{i}_{\mathcal{Z}}(M) & \cong \underset{\fa\in F(\mathcal{Z})}{\varinjlim}H_{\mathfrak{a}}^i(M) \\
 & \cong \underset{\fa\in F(\mathcal{Z})}{\varinjlim}H_{\mathfrak{a}S}^i(M) \\
 & \cong \underset{\fb\in F\left(\widetilde{\mathcal{Z}}\right)}{\varinjlim}H_{\mathfrak{b}}^i(M) \\
 & \cong H^{i}_{\widetilde{\mathcal{Z}}}(M).
\end{split}
\end{equation*}
for every $i \geq 0$.
Hence by part (ii), $N:= H^{i}_{\mathcal{Z}}(M) \cong H^{i}_{\widetilde{\mathcal{Z}}}(M)$ is a $\widetilde{\mathcal{Z}}$-cofinite $S$-module for every $i \geq 0$.
Now, Lemma \ref{3.7} implies that $H^{i}_{\mathcal{Z}}(M)$ is $\mathcal{Z}$-cofinite for every $i\geq 0$.

(iv): Clearly, $\Gamma_{\mathcal{Z}}(M)$ is $\mathcal{Z}$-cofinite.  So by replacing $M$ with $M/\Gamma_{\mathcal{Z}}(M)$,
we can assume that $H^{0}_{\mathcal{Z}}(M)=0$. Now, Lemma \ref{4.2} implies that $\dim\left(\Supp_{R}\left(H^{i}_{\mathcal{Z}}(M)\right)\right) \leq 1$ for all $i\geq 0$, thereby part (iii)
completes the argument.
\end{prf}

In the rest of this section, we apply the technique of way-out functors to prove the main result of this section.

\begin{definition} \label{4.4}
Let $R$ and $S$ be two rings, and $\mathcal{F}: \mathcal{D}(R) \rightarrow \mathcal{D}(S)$ a covariant functor. We say that $\mathcal{F}$ is \textit{way-out left} if for every
$n \in \mathbb{Z}$, there is an $m \in \mathbb{Z}$, such that for any $R$-complex $X$ with $\sup X \leq m$, we have $\sup \mathcal{F}(X) \leq n$.
\end{definition}

The Way-out Lemma appears in \cite[Ch. I, Proposition 7.3]{Ha2}. However, we need a refined version which is tailored to our needs. Let $R$ and $S$ be two rings. Roughly speaking, a triangulated functor is a functor $\mathcal{F}: \mathcal{D}(R) \rightarrow \mathcal{D}(S)$ that preserves shift and distinguished triangles; See \cite[Definition A.7]{CFH}.

\begin{lemma} \label{4.5}
Let $R$ and $S$ be two rings, and $\mathcal{F}: \mathcal{D}(R) \rightarrow \mathcal{D}(S)$ a triangulated covariant functor. Let $\mathcal{A}$ be an additive subcategory of $\mathcal{M}(R)$, and $\mathcal{B}$ an abelian subcategory of $\mathcal{M}(S)$ which is closed under extensions. Suppose that $H_{i}\left(\mathcal{F}(M)\right) \in \mathcal{B}$ for every $M \in \mathcal{A}$ and every $i \in \mathbb{Z}$. If $\mathcal{F}$ is way-out left and $X \in \mathcal{D}_{\sqsubset}(R)$ with $H_{i}(X) \in \mathcal{A}$ for every $i \in \mathbb{Z}$, then $H_{i}\left(\mathcal{F}(X)\right) \in \mathcal{B}$ for every $i \in \mathbb{Z}$.
\end{lemma}

\begin{prf}
See \cite[Lemma 3.2]{DFT}.
\end{prf}

The next result provides us with a suitable transition device from modules to complexes when dealing with cofiniteness.

\begin{lemma} \label{4.6}
Let $\mathcal{Z}$ a stable under specialization subset of $\Spec (R)$. Then the functor ${\bf R}\Gamma_{\mathcal{Z}}(-): \mathcal{D}(R) \rightarrow \mathcal{D}(R)$ is triangulated and way-out left. As a consequence, if $H^{i}_{\mathcal{Z}}(M)$ is $\mathcal{Z}$-cofinite for every finitely generated $R$-module $M$ and every $i \geq 0$, and $\mathcal{M}(R,\mathcal{Z})_{cof}$ is an abelian category, then $H^{i}_{\mathcal{Z}}(X)$ is $\mathcal{Z}$-cofinite for every $X \in \mathcal{D}^{f}_{\sqsubset}(R)$ and every $i \in \mathbb{Z}$.
\end{lemma}

\begin{prf}
It is folklore that if an endofunctor on $\mathcal{D}(R)$ extends from an endofunctor on $\mathcal{M}(R)$, then it commutes with mapping cones. Hence, it can be easily verified that the functor ${\bf R}\Gamma_{\mathcal{Z}}(-): \mathcal{D}(R) \rightarrow \mathcal{D}(R)$ is triangulated. Moreover, if $\sup X \leq n$, then there is a semi-injective resolution $X \xrightarrow{\simeq} I$ of $X$ such that $I_{i}=0$ for every $i \geq n$ due to \cite[Theorem 5.3.26]{CFH}. Applying the functor $\Gamma_{\mathcal{Z}}(-)$ to $I$ and taking homology, we see that $\sup {\bf R}\Gamma_{\mathcal{Z}}(X) \leq n$. It follows that the functor ${\bf R}\Gamma_{\mathcal{Z}}(-)$ is way-out left. Now, let $\mathcal{A}$ be the subcategory of finitely generated $R$-modules, and let $\mathcal{B}:= \mathcal{M}(R,\mathcal{Z})_{cof}$. It can be easily seen that $\mathcal{B}$ is closed under extensions. It now follows from Lemma \ref{4.5} that $H^{i}_{\mathcal{Z}}(X)=H_{-i}\left({\bf R}\Gamma_{\mathcal{Z}}(X)\right) \in \mathcal{B}$ for every $X \in \mathcal{D}^{f}_{\sqsubset}(R)$ and every $i \in \mathbb{Z}$.
\end{prf}

\begin{theorem} \label{4.8}
Let $\mathcal{Z}$ a stable under specialization subset of $\Spec(R)$, and $X \in \mathcal{D}_{\sqsubset}^{f}(R)$. Assume
that either
\begin{enumerate}
\item[(i)]  $R$ is semilocal with $\cd(\mathcal{Z},R) \leq 1$, or
\item[(ii)] $\dim(\mathcal{Z}) \leq 1$, or
\item[(iii)] $\dim\left(\Supp_{R}(X)\right) \leq 2$.
\end{enumerate}
Then $H^{i}_{\mathcal{Z}}(X)$ is $\mathcal{Z}$-cofinite for every $i \in \mathbb{Z}$.
\end{theorem}

\begin{prf}
(i) and (ii): In view of Theorem \ref{3.11} and Lemma \ref{4.3}, one can apply Lemma \ref{4.6} to conclude the claim.

(iii): Let
$$\mathcal{A}:= \left\{M \in \mathcal{M}(R) \suchthat M \text{ is finitely generated and } \dim_{R}(M) \leq 2 \right\},$$
and
$$\mathcal{B}:= \left\{M \in \mathcal{M}(R) \suchthat M \text{ is } \mathcal{Z} \text{-cofinite and } \dim_{R}\left(M\right) \leq 1 \right\}.$$
By Lemma \ref{3.6} and Lemma \ref{4.7} (i), $\mathcal{B}$ is an abelian subcategory of $\mathcal{M}(R)$. In addition, it is closed under extensions.
Now, by Lemmas \ref{4.2} and \ref{4.3} (iv), we have $H^{i}_{\mathcal{Z}}(M) \in \mathcal{B}$ for any $M \in \mathcal{A}$. Considering the triangulated way-out left functor ${\bf R}\Gamma_{\mathcal{Z}}(-)$, Lemma \ref{4.5} yields that $H^{i}_{\mathcal{Z}}(X)=H_{-i}\left({\bf R}\Gamma_{\mathcal{Z}}(X)\right) \in \mathcal{B}$
for every $X \in \mathcal{D}_{\sqsubset}(R)$ with $H_{i}(X) \in \mathcal{A}$ and for every $i \in \mathbb{Z}$.
\end{prf}

\begin{remark} \label{4.9}
Given a stable under specialization subset $\mathcal{Z}$ of $\Spec(R)$, the local cohomology module $H^{i}_{\mathcal{Z}}(X)$ of an $R$-complex $X$ with support in $\mathcal{Z}$, is an all-in-one generalization of the previously known generalized local cohomology modules as outlined in the following discussion.
\begin{enumerate}
\item[(i)] Let $\mathfrak{a}$ be an ideal of $R$, and $M$ and $N$ two $R$-modules. The generalized local
cohomology module $H^{i}_{\mathfrak{a}}(M,N)$ is defined as $$H^{i}_{\mathfrak{a}}(M,N):= \underset{n}{\varinjlim} \Ext^{i}_{R}(M/\mathfrak{a}^{n}M,N)$$ for every $i\geq 0$; see
\cite{Gr} and \cite{He}. It is shown in \cite{Ya} that if $M$ is finitely generated, then we have $H^{i}_{\mathfrak{a}}(M,N)=H^{i}_{\mathcal{Z}}(X)$ for every $i \geq 0$, where $\mathcal{Z}= V(\mathfrak{a})$ and $X={\bf R}\Hom_{R}(M,N)$.
\item[(ii)] Let $\mathfrak{a}$ and $\mathfrak{b}$ be two ideals of $R$, and $M$ an $R$-module.
Let $$W(\mathfrak{a},\mathfrak{b}):=\left\{\mathfrak{p}\in \Spec (R) \suchthat \mathfrak{a}^{n} \subseteq \mathfrak{p}+\mathfrak{b} \text{ for some integer } n\geq 1\right\}.$$
Define a functor $\Gamma_{\mathfrak{a},\mathfrak{b}}(-)$ on $\mathcal{M}(R)$ by setting
$$\Gamma_{\mathfrak{a},\mathfrak{b}}(M):= \left\{x\in M \suchthat \Supp_{R}(Rx) \subseteq W(\mathfrak{a},\mathfrak{b}) \right\},$$
for an $R$-module $M$, and $\Gamma_{\mathfrak{a},\mathfrak{b}}(f):= f|_{\Gamma_{\mathfrak{a},\mathfrak{b}}(M)}$ for an $R$-homomorphism $f:M\rightarrow N$. The generalized local cohomology module $H^{i}_{\mathfrak{a},\mathfrak{b}}(M)$ is defined in \cite{TYY} to be
$H^{i}_{\mathfrak{a},\mathfrak{b}}(M):= R^{i}\Gamma_{\mathfrak{a},\mathfrak{b}}(M)$ for every $i\geq0$.
 It is clear that $H^{i}_{\mathfrak{a},\mathfrak{b}}(M)=H^{i}_{\mathcal{Z}}(X)$ for every $i \geq 0$, where $\mathcal{Z}= W(\mathfrak{a},\mathfrak{b})$ and $X=M$.
\item[(iii)] Let $\Phi$ be a directed poset. By a system of ideals $\varphi$, we mean a family
$\varphi= \{\mathfrak{a}_{\alpha}\}_{\alpha\in\Phi}$ of ideals of $R$, such that $\mathfrak{a}_{\alpha}\subseteq \mathfrak{a}_{\beta}$ whenever $\alpha\geq \beta$, and for any $\alpha,\beta \in \Phi$, there is a $\gamma \in \Phi$ with $\mathfrak{a}_{\gamma} \subseteq \mathfrak{a}_{\alpha}\mathfrak{a}_{\beta}$. Given a system of ideals $\varphi$, define a functor $\Gamma_{\varphi}(-)$ on $\mathcal{M}(R)$ by setting
$$\Gamma_{\varphi}(M):= \left\{x\in M \suchthat \mathfrak{a}x=0  \text{ for some } \mathfrak{a} \in \varphi \right\},$$
for an $R$-module $M$, and $\Gamma_{\varphi}(f):= f|_{\Gamma_{\varphi}(M)}$ for an $R$-homomorphism $f:M\rightarrow N$.
Then the generalized local cohomology module $H^{i}_{\varphi}(M)$ is defined in \cite[Notation 2.2.2]{BS} to be
$H^{i}_{\varphi}(M):=R^{i}\Gamma_{\varphi}(M)$ for every $i\geq0$.
It is easy to see that $H^{i}_{\varphi}(M)=H^{i}_{\mathcal{Z}}(X)$ for every $i \geq 0$, where $\mathcal{Z}=\bigcup_{\mathfrak{a}\in \varphi}V(\mathfrak{a})$ and $X=M$.
\item[(iv)]   Yoshino and Yoshizawa \cite[Theorem 2.10]{YY} have shown that for any abstract local cohomology
functor $\delta:\mathcal{D}_{\sqsubset}(R)
\rightarrow \mathcal{D}_{\sqsubset}(R)$, there is a stable under specialization subset $\mathcal{Z}$ of $\Spec(R)$ such that $\delta\cong {\bf R}\Gamma_{\mathcal{Z}}(-)$.
\end{enumerate}
Accordingly, our Theorems \ref{1.5} and \ref{1.6} generalize the following results:
\begin{enumerate}
\item[(a)] \cite[Proposition 3.6, Corollaries 3.9, 3.10, 3.11, and 3.12]{HV}.
\item[(b)] \cite[Theorem 2.5]{DH}.
\item[(c)] \cite[Theorem 1.3]{DS}.
\item[(d)] \cite[Theorems 1.1 and 1.3]{TGV}.
\item[(e)] \cite[Propositions 6.1 and 7.6, and Corollaries 6.3 and  7.7]{Ha1}.
\item[(f)] \cite[Theorem 2.1]{Ka1}.
\item[(g)] \cite[Theorem 1]{Ka2}.
\end{enumerate}
\end{remark}



\begin{thebibliography}{99}

\bibitem[AB]{AB}{M. Aghapournahr and K. Bahmanpour}, {\it Cofiniteness of general local cohomology modules for small dimensions}, Bull. Korean Math. Soc.,
\textbf{53}(5), (2016), 1341-1352.

\bibitem[AF]{AF}{L. Avramov and H-B. Foxby}, {\it Homological dimensions of unbounded complexes}, J. Pure Appl. Algebra, \textbf{71}(2-3), (1991),
129-155.

\bibitem[BNS]{BNS}{K. Bahmanpour, R. Naghipour and M. Sedghi}, {\it On the category of cofinite modules which is abelian}, Proc. Amer. Math. Soc.,
\textbf{142}(4), (2014), 1101-1107.

\bibitem[BN]{BN}{K. Bahmanpour and R. Naghipour}, {\it Cofiniteness of local cohomology modules for ideals of small dimension}, J. Algebra, \textbf{321}(7),
(2009), 1997-2011.

\bibitem[BS]{BS}{M. Brodmann and R.Y. Sharp}, {\it Local cohomology: An algebraic introduction with geometric applications}, Cambridge Studies in Advanced
Mathematics, \textbf{136}, Cambridge University Press, Cambridge, Second Edition (2013).

\bibitem[Bi]{Bi}{M.H. Bijan-Zadeh}, {\it A common generalization of local cohomology theories}, Glasgow Math. J., \textbf{21}(2), (1980), 173-181.

\bibitem[CFH]{CFH}{L. W. Christensen, H.-B. Foxby, and H. Holm}, {\it Derived category methods in commutative algebra}, draft November 13, 2012.

\bibitem[CW]{CW}{F.C. Cheng and M.Y. Wang}, {\it Homological dimension of G-Matlis dual modules over semilocal rings}, Comm. Algebra, \textbf{21}(4),
(1993), 1215-1220.

\bibitem[DM]{DM}{D. Delfino and T. Marley}, {\it Cofinite modules and local cohomology}, J. Pure Appl. Algebra, \textbf{121}(1), (1997), 45-52.

\bibitem[DFT]{DFT}{K. Divaani-Aazar, H. Faridian, and M. Tousi}, {\it A New Outlook on Cofiniteness}, arXiv:1701.07716 [math.AC].

\bibitem[DH]{DH}{K. Divaani-Aazar and A. Hajikarimi}, {\it Cofiniteness of generalized local cohomology modules for one-dimensional ideals},
Canad. Math. Bull., \textbf{55}(1), (2012), 81-87.

\bibitem[DNT]{DNT}{K. Divaani-Aazar, R. Naghipour and M. Tousi}, {\it Cohomological dimension of certain algebraic varieties}, Proc. Amer.
Math. Soc., \textbf{130}(12), (2002), 3537-3544.

\bibitem[DS]{DS}{K. Divaani-Aazar and R. Sazeedeh}, {\it Cofiniteness of generalized local cohomology modules}, Colloq. Math.,
\textbf{99}(2), (2004), 283-290.

\bibitem[Gr]{Gr}{A. Grothendieck}, {\it Cohomologie locale des faisceaux coh\`{e}rents et th\`{e}or\`{e}mes de Lefschetz locaux et
globaux}, (SGA 2), North-Holland, Amsterdam, 1968.

\bibitem[Ha1]{Ha1}{R. Hartshorne}, {\it Affine duality and cofiniteness}, Invent. Math., \textbf{9}, (1969/1970), 145-164.

\bibitem[Ha2]{Ha2}{R. Hartshorne}, {\it Residues and duality}, Lecture Notes in Mathematics, \textbf{20}, (1966).

\bibitem[HD]{HD}{M. Hatamkhani and K. Divaani-Aazar}, {\it The derived category analogue of the Hartshorne-Lichtenbaum vanishing theorem},
Tokyo. J. Math., \textbf{36}(1), (2013), 195-205.

\bibitem[He]{He}{J. Herzog}, {\it Komplexe, aufl¨osungen und dualit¨at in der lokalen algebra}, Habilitationsschrift, Universit¨at Regensburg, (1970).

\bibitem[HV]{HV}{S. H. Hasanzadeh and A. Vahidi}, {\it On vanishing and cofiniteness of generalized local cohomology modules}, Comm. Algebra,
\textbf{37}(7), (2009), 2290-2299.

\bibitem[HK]{HK}{G. Huneke and J. Koh}, {\it Cofiniteness and vanishing of local cohomology modules}, Math. Proc.
Camb. Phil. Soc., \textbf{110}(3), (1991), 421-429.

\bibitem[HS]{HS}{G. Huneke and R.Y.  Sharp}, {\it Bass numbers of local cohomology modules}, Trans. Amer. Math. Soc.,
\textbf{339}(2), (1993), 765-779.

\bibitem[Ka1]{Ka1}{K. Kawasaki}, {\it On the category of cofinite modules for principal ideals}, Nihonkai Math. J., \textbf{22}(2),
(2011), 67-71.

\bibitem[Ka2]{Ka2}{K. Kawasaki}, {\it Cofiniteness of local cohomology modules for principal ideals}, Bull. London Math. Soc., \textbf{30}(3), (1998), 241-246.

\bibitem[Li]{Li}{J. Lipman}, {\it Lectures on local cohomology and duality}, Local cohomology and its applications, Lecture notes in pure and
applied mathematics, \textbf{226}, (2012), Marcel Dekker, Inc.

\bibitem[Ly]{Ly}{G. Lyubeznik}, {\it Finiteness properties of local cohomology modules (an application of D-modules
to commutative algebra)}, Invent. Math., \textbf{113}(1), (1993), 41-55.

\bibitem[MS]{MS}{A. Mafi and H. Saremi}, {\it On the cofiniteness properties of certain general local cohomology modules}, Acta Sci. Math. (Szeged),
\textbf{74},(3-4), (2008), 501-507.

\bibitem[Me1]{Me1}{L. Melkersson}, {\it Cofiniteness with respect to ideals of dimension one}, J. Algebra, \textbf{372}, (2012), 459-462.

\bibitem[Me2]{Me2}{L. Melkersson}, {\it Modules cofinite with respect to an ideal}, J. Algebra, \textbf{285}(2), (2005), 649-668.

\bibitem[Me3]{Me3}{L. Melkersson}, {\it On asymptotic stability for sets of prime ideals connected with the powers of an ideal}, Math. Proc.
Cambridge Philos. Soc., \textbf{107}(2), (1990), 267-271.

\bibitem[Ro]{Ro}{J.J. Rotman}, {\it An introduction to homological algebra}, Universitext. Springer, New York, second edition, 2009.

\bibitem[Sp]{Sp}{N. Spaltenstein}, {\it Resolutions of unbounded complexes}, Compositio Math., \textbf{65}(2), (1988), 121-154.

\bibitem[TGV]{TGV}{N. Tu Cuong, S. Goto and N. Van Hoang}, {\it On the cofiniteness of generalized local cohomology modules}, Kyoto Journal
of Mathematics, \textbf{55}(1), (2015), 169-185.

\bibitem[TYY]{TYY}{R. Takahashi, Y. Yoshino and T. Yoshizawa}, {\it Local cohomology based on a nonclosed support defined by a pair of ideals},
J. Pure Appl. Algebra, {\bf 213}(4), (2009), 582-600.

\bibitem[WW]{WW}{S. Sather-Wagstaff and R. Wicklein}, {\it Support and adic finiteness for complexes}, Comm. Algebra,
\textbf{45}(6), (2017), 2569-2592.

\bibitem[Ya]{Ya}{S. Yassemi}, {\it Generalized section functors}, J. Pure Appl. Algebra, \textbf{95}(1), (1994), 103-119.

\bibitem[YY]{YY}{Y. Yoshino and T. Yoshizawa}, {\it Abstract local cohomology functors}, Math. J. Okayama Univ., \textbf{53}, (2011), 129-154.

\end{thebibliography}
\end{document}